 \documentclass{amsart}
 \usepackage{euscript}
 \usepackage[top=1in,bottom=1in,left=1.25in,right=1.25in]{geometry} 

 \usepackage{amsfonts}
 \usepackage{amsmath}
 \usepackage{amssymb}
 \usepackage{latexsym,color}
 \usepackage{amsbsy}
 \usepackage{graphicx}
 \usepackage{calc}
 \usepackage{ifthen}
 \usepackage{epsfig}
 \theoremstyle{plain} 
 \newtheorem{thm}{Theorem}[section]
 \newtheorem{con}{Conjecture}[section]
 \newtheorem{lemma}[thm]{Lemma}
 \newtheorem{cor}[thm]{Corollary}
 \newtheorem{pro}[thm]{Proposition}

 \theoremstyle{definition} 
 \newtheorem{example}[thm]{Example}
 \newtheorem{definition}[thm]{Definition}
 \newtheorem{remark}[thm]{Remark}
 \newtheorem{Algorithm}[thm]{Algorithm}

 \newcommand{\ncom}{\newcommand}
 \ncom{\la}{\lambda}
 \ncom{\bm}{\boldmath}
 \ncom{\noi}{\noindent}
 \ncom{\bq}{\begin{equation}}
 \ncom{\eq}{\end{equation}}
 \ncom{\beqn}{\begin{eqnarray*}}
 \ncom{\eeqn}{\end{eqnarray*}}
 \ncom{\ba}{\begin{array}}
 \ncom{\ea}{\end{array}}
 \ncom{\beq}{\begin{eqnarray}}
 \ncom{\eeq}{\end{eqnarray}}
 \ncom{\nno}{\nonumber}
 \ncom{\hs}{\mbox{\hspace{.25cm}}}
 \ncom{\rar}{\rightarrow}
 \ncom{\Rar}{\Rightarrow}
 \ncom{\noin}{\noindent}
 \ncom{\bc}{\begin{center}}
 \ncom{\ec}{\end{center}}
 \ncom{\sz}{\scriptsize}
 \ncom{\fpd}{\Phi(\pi^{'})}
 \ncom{\fp}{\Phi(\pi) }
 \ncom{\nk}{\left< \begin{array}{c}
                        n\\k \end{array} \right>}
 \ncom{\nd}{1^{'},2^{'},\cdots,n^{'}}
 \ncom{\R}{I\!\!R}
 \ncom{\de}{\bigtriangleup (F_{2n}, \leq)}
 \ncom{\del}{\bigtriangleup}
 \ncom{\cov}{<\!\!\!\!\cdot }
 \ncom{\bt}{\begin{thm}}
 \ncom{\bcon}{\begin{con}}
 \ncom{\et}{\end{thm}}
 \ncom{\econ}{\end{con}}
 \ncom{\bl}{\begin{lemma}}
 \ncom{\el}{\end{lemma}}
 \ncom{\bco}{\begin{cor}}
 \ncom{\ds}{\displaystyle}
 \ncom{\eco}{\end{cor}}
 \ncom{\bp}{\begin{pro}}
 \ncom{\ep}{\end{pro}}
 \ncom{\bex}{\begin{example}}
 \ncom{\eex}{\end{example}}
 \ncom{\bd}{\begin{definition}}
 \ncom{\ed}{\end{definition}}
 \ncom{\brm}{\begin{remark}}
 \ncom{\erm}{\end{remark}}
 \ncom{\bal}{\begin{Algorithm}}
 \ncom{\eal}{\end{Algorithm}}

 \ncom{\pf}{\begin{proof}}
 \ncom{\epf}{\end{proof}}
 \ncom{\be}{\begin{enumerate}}
 \ncom{\ee}{\end{enumerate}}
 \ncom{\s}{\subset}
 \ncom{\cc}{\mathcal{C}}
 \ncom{\cf}{\mathcal{F}}
 \ncom{\ac}{{\mathcal{A}}(G,{\mathcal{C}})}
 \ncom{\tc}{{\mathcal{T}}(G,{\mathcal{C}})}
 \ncom{\trc}{{\mathcal{TR}}(G,{\mathcal{C}})}
 \ncom{\kc}{{\mathcal{K}}(G)}
 \ncom{\kcx}{{\mathcal{K}}(G_X')}
 \ncom{\zc}{{\mathcal{Z}}(G)}
 \ncom{\mn}{\mathbb{N}}
 \ncom{\mr}{\mathbb{R}}
 \ncom{\mq}{\mathbb{Q}}
 \ncom{\stack}[2]{\genfrac{}{}{0pt}{}{#1}{#2}} 

 \ncom{\tcr}{\bf\textcolor{red}}
 \ncom{\tcb}{\bf\textcolor{blue}}

 \newenvironment{entry}
   {\begin{list}{}%
      {%
        \setlength{\labelwidth}{35pt}%
        \setlength{\leftmargin}{\labelwidth+\labelsep}%
      }%
   }%
   {\end{list}}

   \newlength{\Mylen}
   \newcommand{\Lentrylabel}[1]{%
     \settowidth{\Mylen}{\emph{#1}}%
     \ifthenelse{\lengthtest{\Mylen > \labelwidth}}%
        {\parbox[b]{\labelwidth}
          {\makebox[0pt][l]{\emph{#1}}\\}}%
        {\emph{#1}}
     \hfil\relax}

 \newenvironment{Lentry}%
   {%
    \begin{entry}}%
   {\end{entry}}

\begin{document}
 \title[Alternating Reachability]{Alternating Reachability}
 \author[Bhattacharya]{Amitava Bhattacharya}
 \address{Bhattacharya: Department of Mathematics, Statistics, and Computer
 Science\\
 University of Illinois at Chicago\\
 Chicago, Illinois 60607-7045, USA\\
 Phone: (312) 413 2163\\
  Fax: (312) 996 1491
 }
 \email{amitava@math.uic.edu}
 \author[Peled]{Uri N. Peled}
 \address{Peled: Department of Mathematics, Statistics, and Computer
 Science\\
 University of Illinois at Chicago\\
 Chicago, Illinois 60607-7045, USA\\
 Phone: (312) 413 2156\\
  Fax: (312) 996 1491}
 \email{uripeled@uic.edu}
 \author[Srinivasan]{Murali K. Srinivasan}
 \address{Srinivasan: Department of Mathematics\\
 Indian Institute of Technology, Bombay\\
 Powai, Mumbai 400076, INDIA\\
 Phone: 91-22-2576 7484\\
  Fax: 91-22-2572 3480}
 \email{mks@math.iitb.ac.in}
 \thanks{UNP and MKS would like to thank Professor Martin Golumbic for his kind
 invitation to visit the Caesarea Edmond Benjamin de Rothschild
 Foundation Institute for Interdisciplinary Applications of Computer
 Science at the University of Haifa, Israel during May--June 2003, where
 part of this work was carried out. The warm hospitality and partial
 support of this visit from CRI is gratefully acknowledged.}

 \dedicatory{Dedicated to the memory of Malka Peled}
 \keywords{colored graphs, alternating walks and trails, Tutte set,
 cycle cone}
 \subjclass[2000]{05C70, 90C27, 90C57}
 \date{October 22, 2005}
 \begin{abstract}
 We consider a graph with colored edges. A trail (vertices may
 repeat but not edges) is called \emph{alternating} when successive
 edges have different colors. Given a set of vertices called
 \emph{terminals}, the \emph{alternating reachability} problem is to
 find an alternating trail connecting distinct terminals, if one exists. A
 special case with two colors is searching for an augmenting path
 with respect to a given matching. In another special case with two
 colors red and blue, the \emph{alternating cone} is defined as the
 set of assignments of nonnegative weights to the edges such that at
 each vertex, the total red weight equals the total blue weight; in
 a companion paper we showed how the search for an integral weight
 vector within a given box in the alternating cone can be reduced to
 the alternating reachability problem in a 2-colored graph. We
 define an obstacle, called a \emph{Tutte set}, to the existence of
 an alternating trail connecting distinct terminals in a colored
 graph, and give a polynomial-time algorithm, generalizing the
 blossom algorithm of Edmonds, that finds either an alternating
 trail connecting distinct terminals or a Tutte set. We use Tutte
 sets to show that an an edge-colored bridgeless graph where each
 vertex has incident edges of at least two different colors has a
 closed alternating trail. A special case with two colors one of
 which forms a matching yields a combinatorial result of Giles and
 Seymour. We show that in a 2-colored graph, the cone generated by
 the characteristic vectors of closed alternating trails is the
 intersection of the alternating cone with the cone generated by the
 characteristic vectors of cycles in the underlying graph.
 \end{abstract}
 \maketitle

 \section{Introduction and Summary}\label{sec1}

 Let $G=(V,E)$ be a graph (we allow parallel edges but not loops).
 A \emph{walk} in $G$ is a sequence
 \beqn
  W &=& (v_0,e_1,v_1,e_2,v_2, \ldots , e_m,v_m),\qquad m \geq 0,
 \eeqn
 where $v_i \in V$ for all $i$, $e_j \in E$ for all $j$, and $e_j$
 has endpoints $v_{j-1}$ and $v_j$ for all $j$. We say that
 $v_1,v_2,\ldots ,v_{m-1}$ are the \emph{internal vertices} of the
 walk $W$. Note that since we are allowing repetitions, the vertices
 $v_0, v_m$ could also be internal vertices. The walk $W$ is said to
 be {\em closed} when $v_0 = v_m$ and is said to be a {\em trail}
 when the edges $e_1,\ldots ,e_m$ are distinct.

 Now assume that the edges of $G$ are colored with a set $C$ of
 colors, where $\#C \geq 2$, the coloring being given by ${\mathcal
 C}: E \rar C$. We say that $(G,\mathcal{C})$ is an
 \emph{edge-colored graph}. The walk $W$ above is {\em internally
 alternating} when $\mathcal{C}(e_j) \neq \mathcal{C}(e_{j+1})$ for
 each $j = 1,\ldots, m-1$, and is {\em alternating} when in addition
 if $W$ is closed then $\mathcal{C}(e_m) \neq \mathcal{C}(e_1)$
 (note that a walk can be closed and internally alternating without
 being alternating, but if $v_0 \neq v_m$, there is no distinction
 between internally alternating and alternating walks and we use the
 word alternating in this case). A closed alternating walk
 (respectively, trail) is abbreviated as \emph{CAW} (respectively,
 \emph{CAT}).

 Given an edge-colored graph $G=(V,E)$ and a set $S \subseteq V$ of
 vertices called \emph{terminals}, the \emph{alternating
 reachability problem} is to either find an alternating trail
 connecting distinct terminals or show that none exists. Our
 motivation for considering this problem arises from the following
 two special cases with two colors, say red and blue.

 \textbf{(a)} Let $G=(V,E)$ be a simple graph and let $M \subseteq
 E$ be a matching. Color the edges in $M$ red and the edges in $E -
 M$ blue, and let $S$ be the set of exposed vertices of $M$. The
 alternating reachability problem for these data is equivalent to
 finding an augmenting path with respect to $M$, if one exists.

 \textbf{(b)} Let $G=(V,E)$ be a graph whose edges have been colored
 red and blue. The {\em alternating cone} is defined as the set of
 assignments of nonnegative weights to the edges such that at each
 vertex, the total red weight equals the total blue weight (we
 define the alternating cone more formally below). In a companion
 paper \cite{BPS1} we showed how the search for an integral weight
 vector within a given box in the alternating cone can be reduced to
 the alternating reachability problem in a 2-colored graph.

 We now outline the main results of this paper.

 In Section~\ref{sec4} we show, generalizing the approach of Edmonds
 \cite{e} as explained in Lov\'{a}sz and Plummer \cite{lp}, that
 given an edge-colored graph $G=(V,E)$ and a set of terminals,
 either there is an alternating trail connecting distinct terminals,
 or else there is a subset of nonterminals, called a \emph{Tutte
 set}, that acts as an obstruction to such alternating trails (for a
 precise definition of a Tutte set, see Section~\ref{sec4}).
 Moreover, we present a polynomial-time algorithm that finds one or
 the other. The alternating reachability problem (in a slightly
 different version) was first  considered by Tutte, in a
 nonalgorithmic form, in his book on graph theory \cite{t1} and in
 the paper \cite{t2}. Tutte \cite{t2} called the obstructions to
 alternating trails connecting distinct terminals $r$-barriers.
 Every $r$-barrier is a Tutte set but not conversely. It turns out
 that there is a very minor error in Tutte's work: Tutte's theory
 actually produces Tutte sets in our sense and not $r$-barriers in
 his sense. Also, it is easy to give an example where there exists a
 unique Tutte set, which is not an $r$-barrier.

 In Section~\ref{sec5} we use Tutte sets to prove a combinatorial
 result on CAT's in edge-colored bridgeless graphs. To motivate this
 result consider the following statements:
 \be
 \item[(i)] Let $D$ be a directed graph in which every vertex has positive
 indegree and outdegree. Then $D$ has a directed circuit (this is easily
 proved).
 \item[(ii)] Let $(G,\mathcal{C})$ be an edge-colored graph such that for
 every vertex $v$ of $G$, we can find edges of two different colors incident
 with $v$.  Then $(G,\mathcal{C})$ has a CAW (this can be proved by a simple
 alternating walk argument similar to the proof of Theorem~2.2 of \cite{BPS1}).
 \ee
 Theorem~\ref{catbg} strengthens the hypothesis of statement (ii)
 above by assuming in addition that $G$ is bridgeless. It concludes,
 using Tutte sets, that $(G,\mathcal{C})$ has a CAT. We also deduce
 a result (Theorem~\ref{Seymour-Giles}) due to Giles and Seymour
 \cite{s} on cycles in bridgeless graphs from a special case of
 Theorem~\ref{catbg}, where there are two colors, red and blue, and
 the red edges form a matching.

 In Section~\ref{sec6} we use Theorem~\ref{catbg} to prove an intersection
 theorem for certain polyhedral cones associated to {\em 2-colored graphs}
 (i.e., edge-colored graphs where the number of colors is two) that were
 defined in \cite{BPS1}. We now recall some definitions from \cite{BPS1}.

 Let $G=(V,E)$ be a graph.
 Assume that the edges of $G$ are colored red or blue, the coloring
 being given by ${\mathcal C}: E \rar \{R,B\}$.
 Consider the real
 vector space ${\mathbb R}^E$, with coordinates indexed by the set
 of edges of $G$. We write an element $x \in {\mathbb R}^E$ as $x=
 (x(e) : e \in E)$. For an edge $e \in E$, the characteristic vector
 $\chi(e) \in {\mathbb R}^E$ is defined by
 $\chi(e)(f) = \left\{ \ba{cc} 1, & \mbox{if }f=e \\
                                          0, & \mbox{if }f \neq e.
                                   \ea \right.$

\begin{figure}
 \centerline{
 \scalebox{.6}
 {
 \epsfclipon
 \epsffile[249 606 408 752]{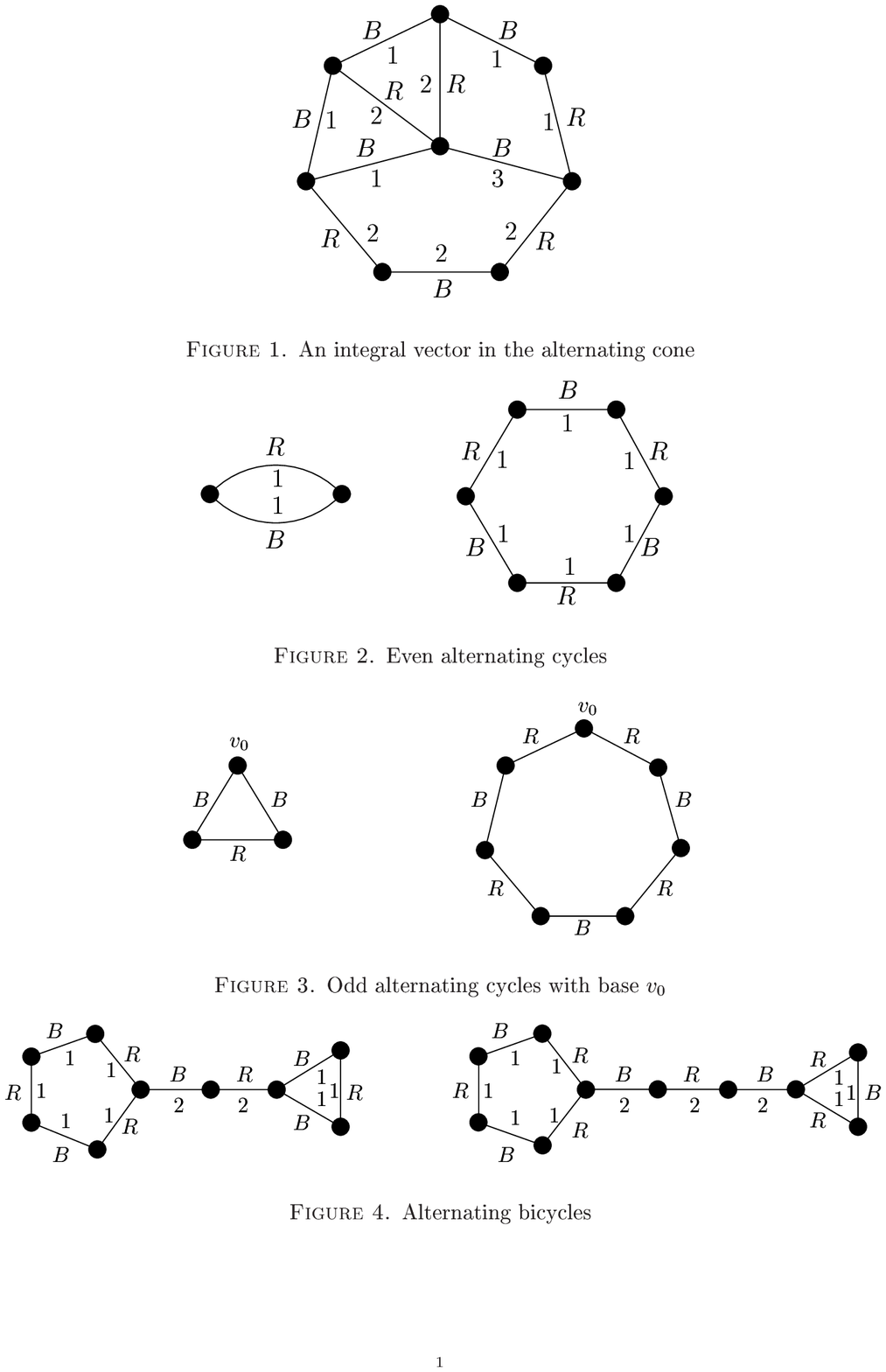}} }
 \caption{An integral
 vector in the alternating cone} \label{fig1_1}
 \end{figure}

 The \emph{cone of closed alternating walks} or simply the
 \emph{alternating cone} ${\mathcal A}(G,{\mathcal C})$ of a
 2-colored graph $(G,\mathcal{C})$ (when the coloring ${\mathcal C}$ is
 understood, we suppress it from the notation and write ${\mathcal
 A}(G)$), is defined to be the set of all vectors $x = (x(e) : e \in
 E )$ in ${\mathbb R}^E$ satisfying the following system of
 homogeneous linear inequalities:
 \beq \label{3}
 \sum_{e \in E_R(v)} x(e) - \sum_{e \in E_B(v)} x(e) & = & 0, \qquad v \in V, \\
 \label{4}                                           x(e) & \geq &
 0,\qquad e \in E .
 \eeq
 We refer to (\ref{3}) as the \emph{balance condition} at vertex
 $v$.
 Figure~\ref{fig1_1} illustrates a 2-colored graph together
 with an integral vector in its alternating cone.
 In~\cite{BPS1} we determined the extreme rays and dimension of the
 alternating cone and showed that searching for an
 integral vector within a given box in the alternating cone can be
 reduced to searching for an alternating trail connecting two given
 vertices in a residual 2-colored graph.

\begin{figure}
  \begin{center}
  \scalebox{.6}
  {\epsfclipon
  \epsffile[90 605 470 760]{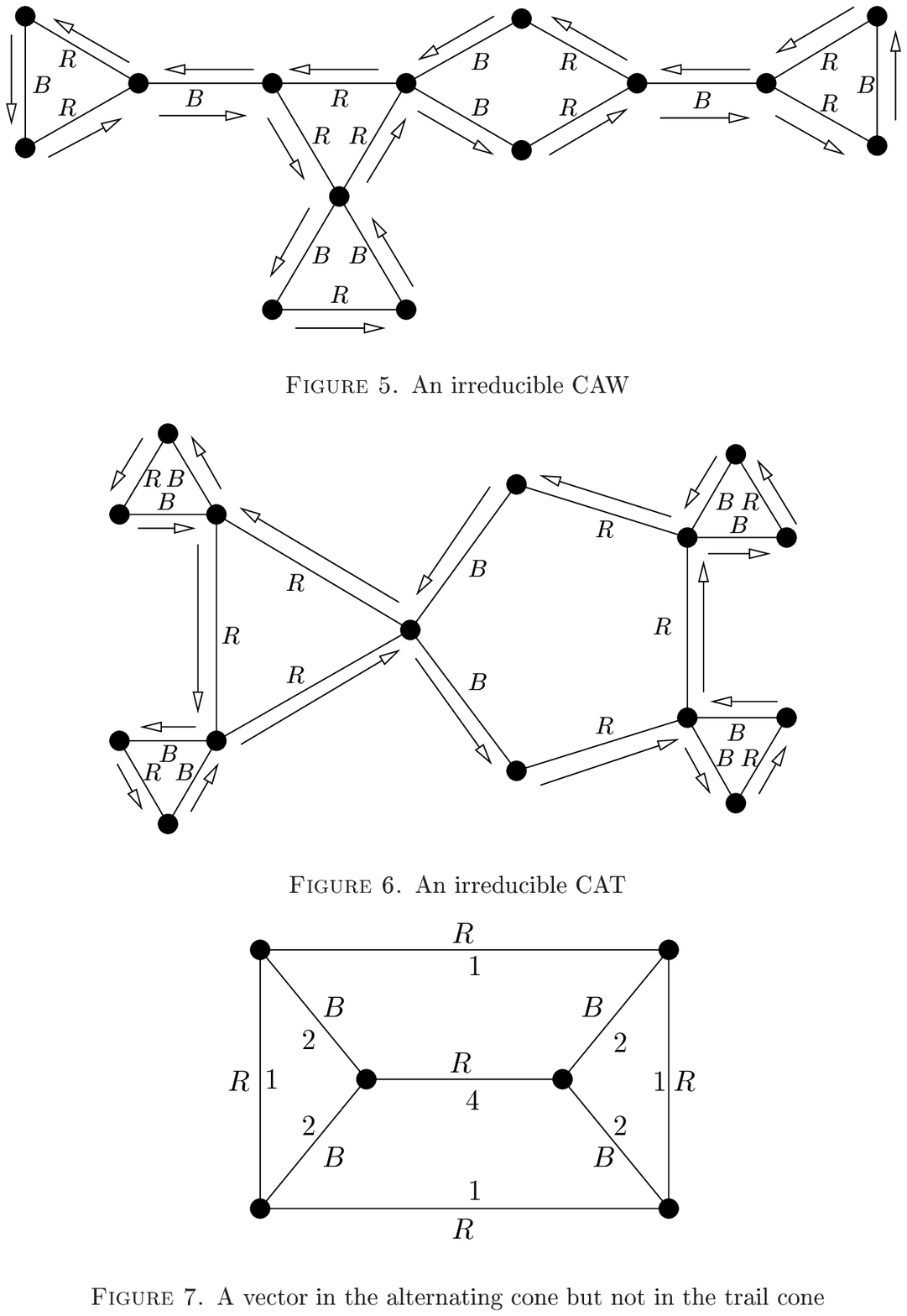}}
   \caption{An irreducible CAW}
   \label{fig1_5}
  \end{center}
 \end{figure}

\begin{figure}
  \begin{center}
  \scalebox{.6}
  {\epsfclipon
  \epsffile[140 395 435 585]{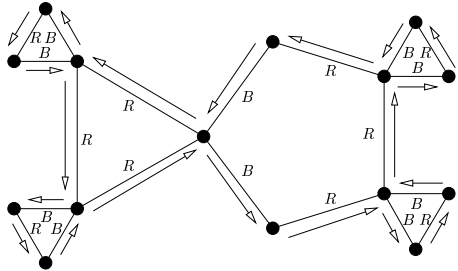}}
   \caption{An irreducible CAT}
   \label{fig1_6}
  \end{center}
 \end{figure}

 The characteristic vector of a walk $W = (v_0,e_1,v_1,e_2,v_2,
 \ldots , e_m,v_m)$ is defined to be $\chi(W) = \sum_{i=1}^m
 \chi(e_i)$. A CAW $W$ in a 2-colored graph is said to be
 \emph{irreducible} if $\chi(W)$ cannot be written as $\chi(W_1) +
 \chi(W_2)$ for any CAW's $W_1$ and $W_2$. Similarly, a A CAT $T$ is
 said to be \emph{irreducible} if $\chi(T)$ cannot be written as
 $\chi(T_1) + \chi(T_2)$ for any CAT's $T_1$ and $T_2$.
 Figure~\ref{fig1_5} depicts an irreducible CAW (with direction of
 walk indicated by an arrow) and Figure~\ref{fig1_6} depicts an
 irreducible CAT. Irreducibility is easily seen. A simple
 alternating walk argument (see \cite{BPS1}) shows that every
 integral vector in the alternating cone is a nonnegative integral
 combination of characteristic vectors of irreducible CAW's.


\begin{figure}
  \begin{center}
  \scalebox{.6}
  {\epsfclipon
  \epsffile[180 235 395 380]{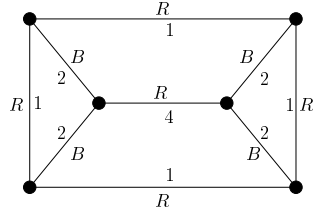}}
   \caption{A vector in the alternating cone but not in the trail cone}
   \label{fig1_7}
  \end{center}
 \end{figure}

\begin{figure}
  \begin{center}
  \scalebox{.6}
  {\epsfclipon
  \epsffile[45 600 597 765]{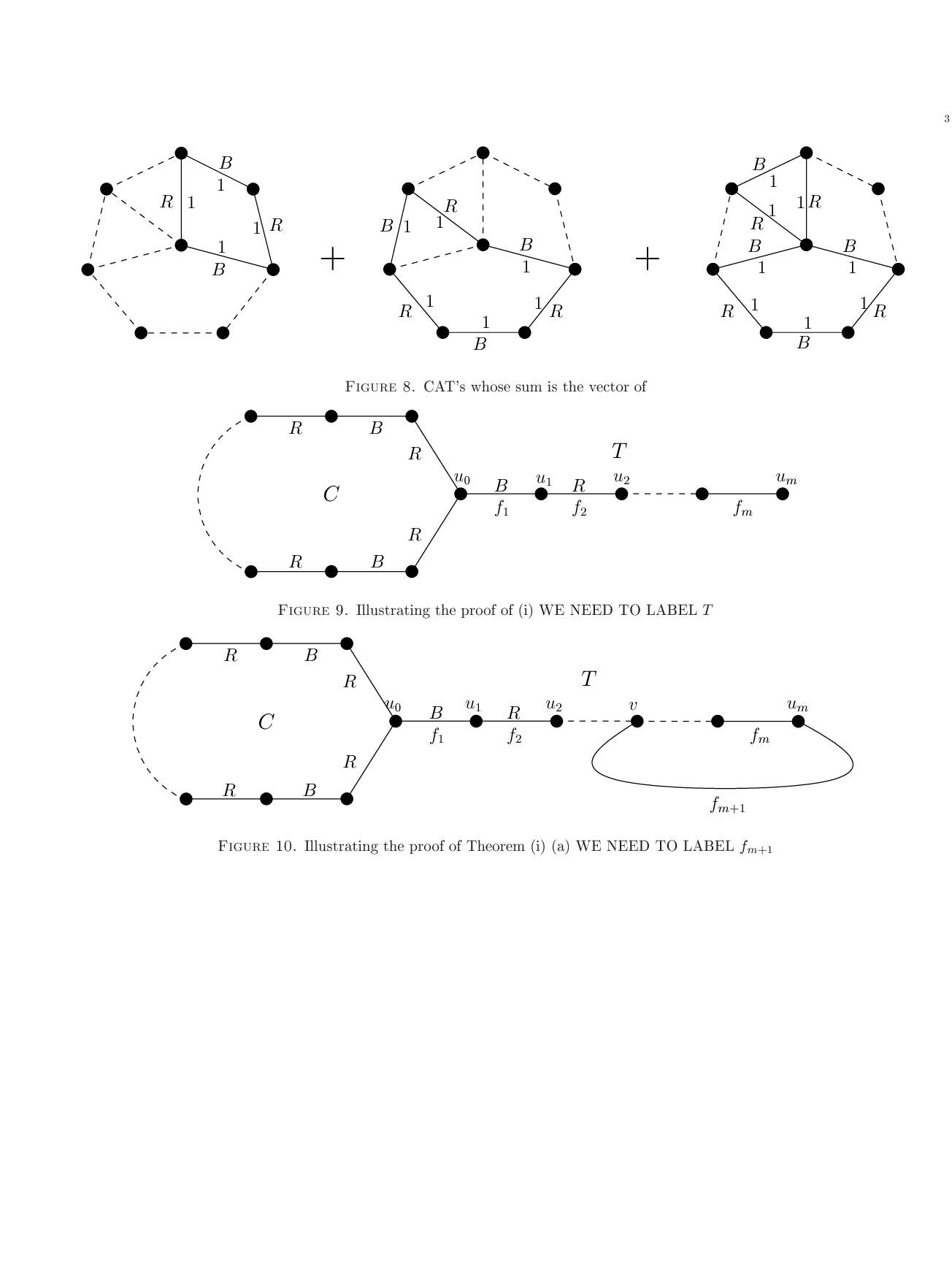}}
   \caption{CAT's whose sum is the vector of Figure \protect\ref{fig1_1}}
   \label{fig1_8}
  \end{center}
 \end{figure}

 Circulations in directed graphs can be thought of in terms of flows
 along the arcs obeying the conservation constraint at every vertex.
 For example, the characteristic vector of a directed circuit
 corresponds to a unit of flow along the circuit. Such an
 interpretation is not available in the case of vectors in the
 alternating cone. The irreducible CAW of Figure~\ref{fig1_5} does
 not correspond to a flow  in an intuitive sense. On the other hand,
 the characteristic vector of an irreducible CAT \emph{can} be
 thought of as a unit of flow around the trail. For a 2-colored
 graph $G=(V,E),\; {\mathcal C}: E\rar \{R,B\}$, it is thus natural
 to consider the convex polyhedral cone ${\mathcal T}(G,{\mathcal
 C}) \subseteq {\mathbb R}^E$ generated by  the characteristic
 vectors of the CAT's in $(G,\mathcal{C})$. We call ${\mathcal
 T}(G,{\mathcal C})$ the \emph{cone of closed alternating trails},
 or simply the \emph{trail cone}, of $(G,\mathcal{C})$. For example,
 it is easily seen that the integral vector in the alternating cone
 from Figure~\ref{fig1_7} is not in the trail cone. On the other
 hand, the integral vector in the alternating cone from
 Figure~\ref{fig1_1} can be written as a sum of characteristic
 vectors of three CAT's, as shown in Figure~\ref{fig1_8}.

 Consider a CAT in a 2-colored graph. Its characteristic vector
 satisfies the balance condition at every vertex. If we ignore the
 colors, the edge set of the CAT is a disjoint union of the edge
 sets of some cycles in the underlying graph. This shows that a
 nonnegative integral combination (that is to say, a linear
 combination with nonnegative integral \emph{coefficients}) of
 characteristic vectors of CAT's satisfies the balance condition at
 every vertex and can be written as a nonnegative integral
 combination of characteristic vectors of cycles in the underlying
 graph. We conjecture that the converse of this observation is also
 true:
 \bcon \label{conj}
 Let $G=(V,E),\; {\mathcal C}: E\rar \{R,B\}$ be a 2-colored graph
 and let $y \in {\mathbb N}^E$.
 \noi Then $y$ is a nonnegative integral combination of
 characteristic vectors of CAT's in $(G,\mathcal{C})$ if and only if
 \be
 \item[\emph{(i)}] $y$ satisfies the balance condition at every vertex, i.e.,
 $y \in {\mathcal A}(G, {\mathcal C})$;

 \item[\emph{(ii)}] $y$ can be written as a nonnegative integral combination
 of characteristic vectors of cycles in $G$.
 \ee
 \econ
 Let ${\mathcal Z}(G)$ denote the cone in ${\mathbb R}^E$ generated
 by the characteristic vectors of the cycles in $G$. Seymour
 \cite{s} found  the linear inequalities determining ${\mathcal
 Z}(G)$ (Theorem~\ref{sey}). The observation in the paragraph
 preceding Conjecture~\ref{conj} shows that ${\mathcal
 T}(G,{\mathcal C}) \subseteq {\mathcal A}(G,{\mathcal C}) \cap
 {\mathcal Z}(G)$. From Conjecture~\ref{conj} it is easy to show
 that ${\mathcal T}(G,{\mathcal C}) = {\mathcal A}(G,{\mathcal C})
 \cap {\mathcal Z}(G)$ : take a rational vector $y  \in {\mathcal
 A}(G,{\mathcal C}) \cap {\mathcal Z}(G)$; for a suitably large
 positive integer $k$, $ky$ satisfies conditions (i) and (ii) of the
 conjecture; by the conjecture $ky$ (respectively, $y$) is a
 nonnegative integral (respectively, nonnegative rational)
 combination of characteristic vectors of CAT's. In
 Section~\ref{sec6} we prove that indeed ${\mathcal T}(G,{\mathcal
 C}) = {\mathcal A}(G,{\mathcal C}) \cap {\mathcal Z}(G)$. Our proof
 is an adaptation of Seymour's argument for ${\mathcal Z}(G)$.
 Seymour's inductive proof is based on the Giles-Seymour lemma
 (Theorem~\ref{Seymour-Giles}) and  likewise, our
 inductive proof uses Theorem~\ref{catbg}.

 We remark that in this paper we focus on graph-theoretical aspects
 of the alternating cone and not on algorithmic efficiency. We do
 consider algorithms, but always with a view to obtaining
 graph-theoretical results.

 We now collect in one place certain commonly used definitions in
 the rest of the paper. Let $G=(V,E)$ be a graph and consider a walk
 \beq
 \label{w}
 W &=& (v_0,e_1,v_1,e_2,v_2, \ldots , e_m,v_m),\qquad m \geq 0,
 \eeq
 in $G$. We say that $W$ is a
 $v_0$-$v_m$ walk of length $m$. We call $e_1$ the \emph{first} edge
 of $W$ and $e_m$ the \emph{last} edge of $W$.
 The walk $W^R$ is the
 $v_m$-$v_0$ walk obtained by reversing the sequence (\ref{w}).
 The walk $W$ is said to be
 \begin{Lentry}
 \item[a path] when the edges $e_1,\ldots ,e_m$ are distinct and the
 vertices $v_0,\ldots ,v_m$ are distinct;

 \item[a cycle] when $W$ is closed, the edges $e_1,\ldots ,e_m$ are
 distinct, and the vertices $v_0,\ldots ,v_{m-1}$ are distinct.
 \end{Lentry}
 We have defined paths and cycles as special classes of walks.
 However, sometimes it is more convenient to  think of paths and
 cycles as subgraphs, as is done usually. This will be clear from
 the context. If $W_1$ is a $u$-$v$ walk and $W_2$ is a $v$-$w$
 walk, then the \emph{concatenation} of $W_1$ and $W_2$, denoted
 $W_1*W_2$, is the $u$-$w$ walk obtained by walking from $u$
 to $v$ along $W_1$ and continuing by walking from $v$ to $w$
 along $W_2$. Note that if $W_1$ and $W_2$ are trails, then
 $W_1*W_2$ is a trail whenever $W_1$ and $W_2$ have no edges in
 common.

 Now let $(G,\mathcal{C})$ be a 2-colored graph. The walk $W$ in
 (\ref{w}) is said to be
 \begin{Lentry}
 \item[an even alternating cycle] when $W$ is a cycle of even length and
 $W$ is alternating; an even alternating cycle will also be called
 simply an \emph{alternating cycle};

 \item[an odd internally alternating cycle with base $v_0$] when $W$ is
 a $v_0$-$v_0$ cycle of odd length and $W$ is internally
 alternating.
 \end{Lentry}

 \section{Alternating Trails in an Edge-Colored Graph}\label{sec4}

 Let $G=(V,E)$ be a graph and let
 $\cc : E \rightarrow C$ be an edge coloring.
 In this section we
 consider the \emph{alternating reachability problem}: given a set
 $S$ of vertices called \emph{terminals}, either find an alternating
 trail connecting distinct terminals or show that no such trail
 exists. For the rest of this section we consider $G$, ${\mathcal
 C}$ and $S$ as fixed.

 The problem of finding a CAT through a given edge $e$ in an
 edge-colored graph can be easily reduced to the alternating
 reachability problem: let $e = \{s,t\}$, $s \neq t$. Remove $e$
 from the graph, add two new vertices $s'$ and $t'$, add a new edge
 with the color $\mathcal{C}(e)$ between  $s'$ and $s$ and one
 between $t'$ and $t$, and let $S = \{s',t'\}$. Clearly the
 alternating reachability problem in the new graph is equivalent to
 the original problem.

 As mentioned in the introduction, the problem of finding an
 augmenting path with respect to a given matching is also reducible
 to the alternating reachability problem.

 The alternating reachability problem was first considered by Tutte
 \cite{t1,t2}. We discuss Tutte's work at the end of this
 section. Our solution to the alternating reachability problem is
 along the lines of the blossom forest algorithm of Edmonds \cite{e},
 as explained in Section~9.1 of Lov\'{a}sz and Plummer's book
 \cite{lp}. The solution is in terms of Tutte sets, defined below.
 \bd
 \label{dfs}
 A subset $A \subseteq (V - S)$ is a \emph{Tutte set} when
 \be
 \item [(i)] each component of $G - A$ has at most one terminal;
 \item [(ii)] $A$ can be written as a
 disjoint union (denoted $\dot{\cup}$, empty blocks allowed)
 \[A = \dot{\bigcup}_{c \in C} A(c)\]
 such that conditions (a), (b), and (c) below hold.

 \noi A vertex $u \in A$ is said to have \emph{color} $c$ if $u \in A(c)$
 (there is a unique such $c$). An edge $e \in E$ is said to be
 \emph{mismatched} if $e$ connects a vertex $u \in A$ with a vertex
 $v \in V - A$ and $\cc(e)$ is different from the color of $u$, or
 $e$ connects two vertices $u,v \in A$ and $\cc(e)$ is different
 from both the colors of $u$ and of $v$.

 \noi Conditions (a), (b), and (c) are as follows:
 \be
 \item [(a)] if $H$ is a component of $G-A$ containing a terminal,
 then there is no mismatched edge with an endpoint in $H$;
 \item [(b)] if $H$ is a component of $G-A$ containing no terminals,
 then there is at most one mismatched  edge with an endpoint in $H$;
 \item [(c)] there are no mismatched edges with both endpoints in $A$.
 \ee
 \ee
 \ed

 The next theorem shows that a Tutte set is an obstruction to the existence
 of an alternating trail connecting distinct terminals.
 \bt \label{ts} Suppose a Tutte set $A$ exists. Let $s$ be a
 terminal, let $H$ be a component of $G-A$ containing the vertex $t$
 but not $s$, and assume that there is no mismatched edge with an
 endpoint in $H$. Then there is no alternating $s$-$t$ trail. In
 particular, there is no alternating trail connecting distinct
 terminals. \et
 \pf
 We assert that if an alternating trail $T$ enters $A$ from $V-A$ via an edge
 that is not mismatched, then the next time $T$ leaves $A$, it can
 only be via a mismatched edge. Specifically, suppose that
 \[T = (v_0,e_1,v_1,e_2,v_2, \ldots , e_m,v_m)\]
 is an alternating trail with $v_0  \in V-A$, $v_1  \in A$, and
 ${\mathcal C}(e_1) = \text{ color of } v_1$. Assume that there
 exists some $j \geq 2$ with $v_j \in V-A$, and let $i$ be the least
 such $j$. Proving the assertion amounts to showing
 $\mathcal{C}(e_i) \neq \text{ color of } v_{i-1}$.

 We show by induction that all $l=1,\ldots i-1$ satisfy $\cc(e_l) =
 \text{ color of } v_l$. The base case $l = 1$ follows from
 hypothesis. Now assume that the statement holds for $l = t$, where
 $t< i-1$. Thus $\cc(e_t) = c$, where $c$ is the color of $v_t$.
 Since  $T$ is alternating, $d = \cc(e_{t+1}) \neq c$. Since
 $v_{t+1} \in A$, it follows from condition (ii)(c) in
 Definition~\ref{dfs} that the color of $v_{t+1}$ is $d$.

 We have shown that $\cc(e_{i-1}) = \text{ color of } v_{i-1}$.
 Since $T$ is alternating, $\cc(e_i) \neq \cc(e_{i-1})$, and thus
 $\mathcal{C}(e_i) \neq \text{ color of } v_{i-1}$, which proves the
 assertion.

 Now suppose that $T$ is an alternating $s$-$t$ trail. Since $s$ and
 $t$ are in different components of $G-A$, $T$ must enter $A$. By
 (ii)(a) in Definition~\ref{dfs}, the first time $T$ enters $A$, it
 must be via an edge that is not mismatched. Since $t \notin A$, $T$
 must leave $A$. By the assertion, the first time $T$ leaves $A$, it
 is via a mismatched edge. By hypothesis, the component $K$ of $G-A$
 that $T$ enters upon leaving $A$ (for the first time) cannot be the
 destination component $H$ containing $t$, and thus $T$ must leave
 $K$ and enter $A$ again. By (ii)(b) in Definition~\ref{dfs}, this
 entry must be via an edge that is not mismatched (since the only
 mismatched edge has already been used for entering $K$). Therefore,
 by the assertion, when $T$ leaves $A$ for the next time, it must be
 via a mismatched edge. Continuing this argument we see that every
 time $T$ leaves $A$, it must be via a mismatched edge. Thus $T$ can
 never reach the destination component $H$ containing $t$, a
 contradiction. \epf

 For later use we record the following lemma.
 \bl
 \label{os}
 Suppose a Tutte set $A$ exists. If $s$ is a terminal and $u \in
 A(c)$, then the last edge of each alternating $s$-$u$ trail has
 color $c$. \el
 \pf
 Let
 \[T = (v_0,e_1,v_1,e_2,v_2, \ldots , e_m,v_m)\]
 be an alternating trail with $s=v_0$ and $u=v_m$. By (ii)(a) in
 Definition~\ref{dfs}, we see that when $T$ enters $A$ for the first
 time, it must be via an edge that is not mismatched. Just as in
 Lemma~\ref{ts}, we can now show that every time $T$ leaves $A$, it
 must be via a mismatched edge, and every time $T$ enters $A$, it
 must be via a edge that is not mismatched. Since $u \in A$, $T$
 must enter $A$ for the last time, say via the edge $e_l$. Since
 $e_l$ is not mismatched, $\cc(e_l) = \mbox{ color of } v_l$. The
 induction argument used in proving the assertion in the proof of
 Theorem~\ref{ts} now shows that $\cc(e_m) = \text{ color of } v_m =
 c$.
 \epf

 We now want to show that the converse of Theorem~\ref{ts} holds,
 i.e., if there is no alternating trail connecting distinct
 terminals, then a Tutte set exists. We shall present an algorithm
 that finds such a trail or a Tutte set. First, we need to make a
 few definitions.

 For a non root vertex $u$ in a rooted forest, the \emph{predecessor
 edge} of $u$ is the first edge in the unique path from $u$ to the
 root of the component containing $u$. Given a partition $\pi$ of
 the vertex set of $G$, by the \emph{shrunken graph} $G \times \pi$
 we mean the graph obtained from $G$ by shrinking each block of
 $\pi$ into a vertex and discarding loops. In other words, the
 vertex set of $G\times \pi$ is the set of blocks of $\pi$, the edge
 set of $G\times \pi$ is the set of edges of $G$ whose endpoints lie
 in different blocks of $\pi$, and the endpoints of an edge $e$ (in
 $G\times \pi$) are the blocks of $\pi$ in which the endpoints of
 $e$ (in $G$) lie. For $v \in V$, the block of $\pi$ containing $v$
 is denoted by $[v]$. For a subset $U \subseteq V$, the subgraph of
 $G$ induced on $U$ is denoted by $G[U]$.

 \bd
 \label{db}
 A subgraph $G' = (V', E')$ of $G$ is said to
 be a \emph{blossom with base $u$} when
 \be
 \item [(i)] $u \in V'$;
 \item [(ii)] for each $v \in V'$ such that $v \neq u$, $G'$ has two $v$-$u$
 alternating trails whose first edges have different colors.
 \ee
 \ed

 Condition (ii) allows us to extend any alternating trail reaching $G'$ up
 to $u$.

 Note that a subgraph consisting of a single vertex $u$ is a blossom with
 base $u$.

 \bd \label{dcb} Let $c \in C$. A subgraph $G' = (V', E')$ of $G$ is
 said to be a \emph{$c$-blossom with base $u$} when
 \be
 \item [(i)] $u \in V'$;
 \item [(ii)] for each $v \in V'$ such that $v \neq u$, $G'$ has
 two $v$-$u$ alternating trails whose first edges have different
 colors and whose last edges have colors different from $c$;
 \item [(iii)] $G'$ has a $u$-$u$ internally alternating trail
 of positive length whose first and last edges have
 colors different from $c$.
 \ee
 \ed

 Conditions (ii) and (iii) allow us to extend any alternating trail
 reaching $G'$ (even a trail reaching $u$ via an edge of color $c$)
 up to $u$ and then be ready to continue with color $c$.

\begin{figure}
  \begin{center}
  \scalebox{.6}
  {\epsfclipon
  \epsffile[265 590 405 750]{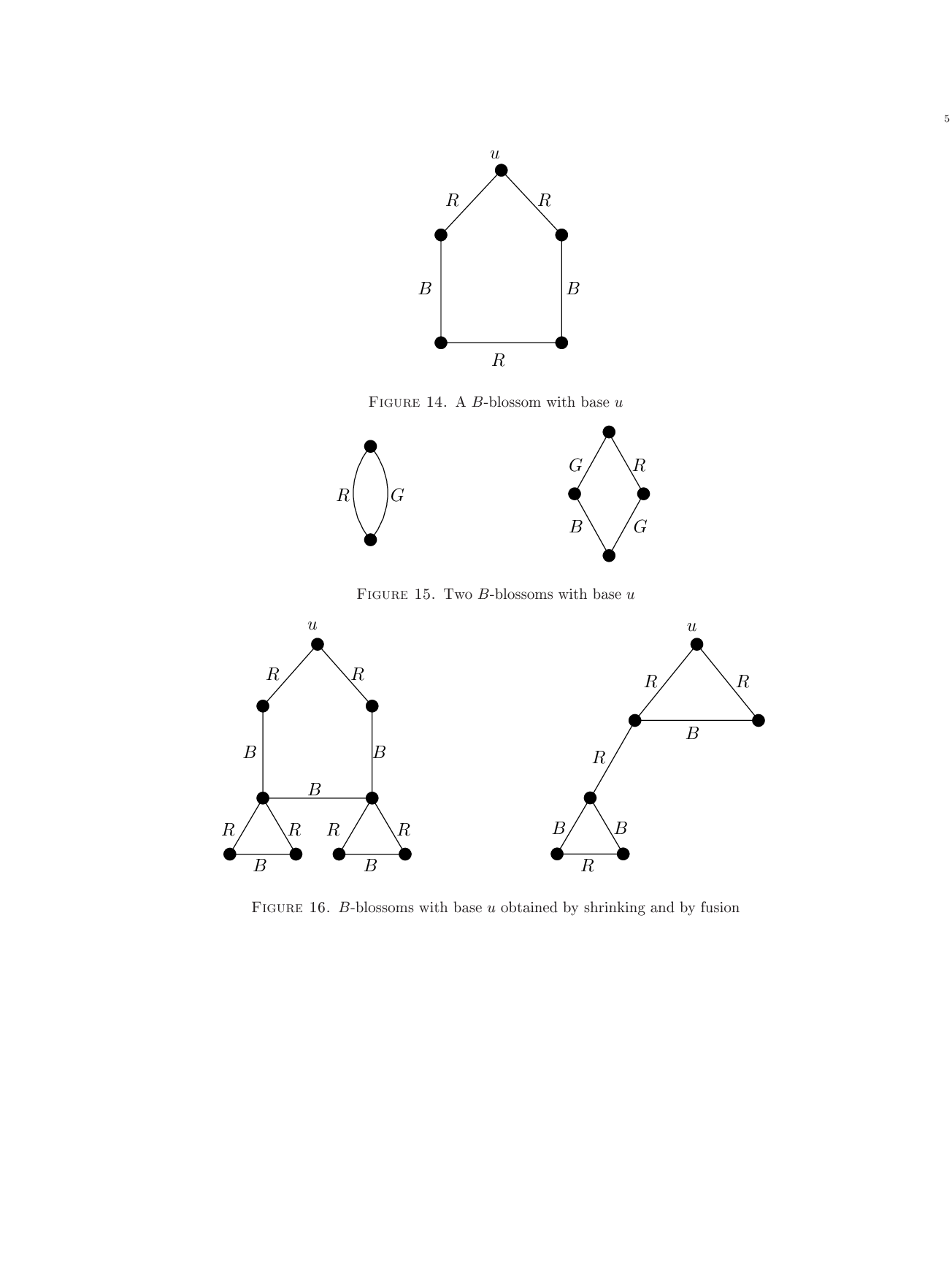}}
   \caption{A $B$-blossom with base $u$}
   \label{fig_15}
  \end{center}
 \end{figure}

 Note that for all $c$, a subgraph consisting of a single vertex $u$
 is \emph{not} a $c$-blossom.
 \bex
 The following are examples of blossoms.
 \be
 \item [(i)] If there are only two colors $R$ and $B$, an odd internally alternating
 $u$-$u$ cycle with two $R$ edges incident at $u$ is a $B$-blossom
 (see Figure~\ref{fig_15}).

\begin{figure}
  \begin{center}
  \scalebox{.6}
  {\epsfclipon
  \epsffile[215 465 440 570]{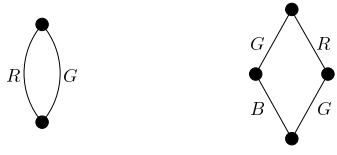}}
   \caption{Two $B$-blossoms with base $u$}
   \label{fig_16}
  \end{center}
 \end{figure}

 \item [(ii)] Assuming there are three colors $R$, $B$ and $G$,
 Figure~\ref{fig_16} depicts two $B$-blossoms with base $u$.

\begin{figure}
  \begin{center}
  \scalebox{.6}
  {\epsfclipon
  \epsffile[135 259 515 440]{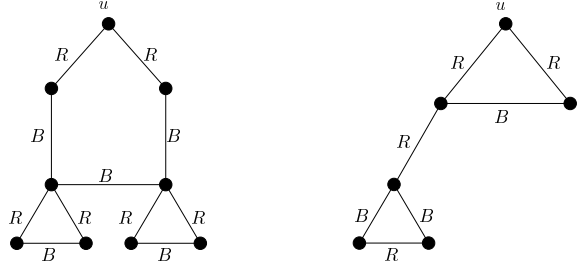}}
   \caption{$B$-blossoms with base $u$ obtained by shrinking and by fusion}
   \label{fig_17}
  \end{center}
 \end{figure}

 \item [(iii)] Figure~\ref{fig_17} depicts two more $B$-blossoms with base $u$. At this
 point this can be verified directly from Definition~\ref{dcb}. Later we
 shall see that the first of these $B$-blossoms arises by shrinking and
 the second by fusion.

 \ee
 \eex

\begin{figure}
  \begin{center}
  \scalebox{.6}
  {\epsfclipon
  \epsffile[85 480 500 760]{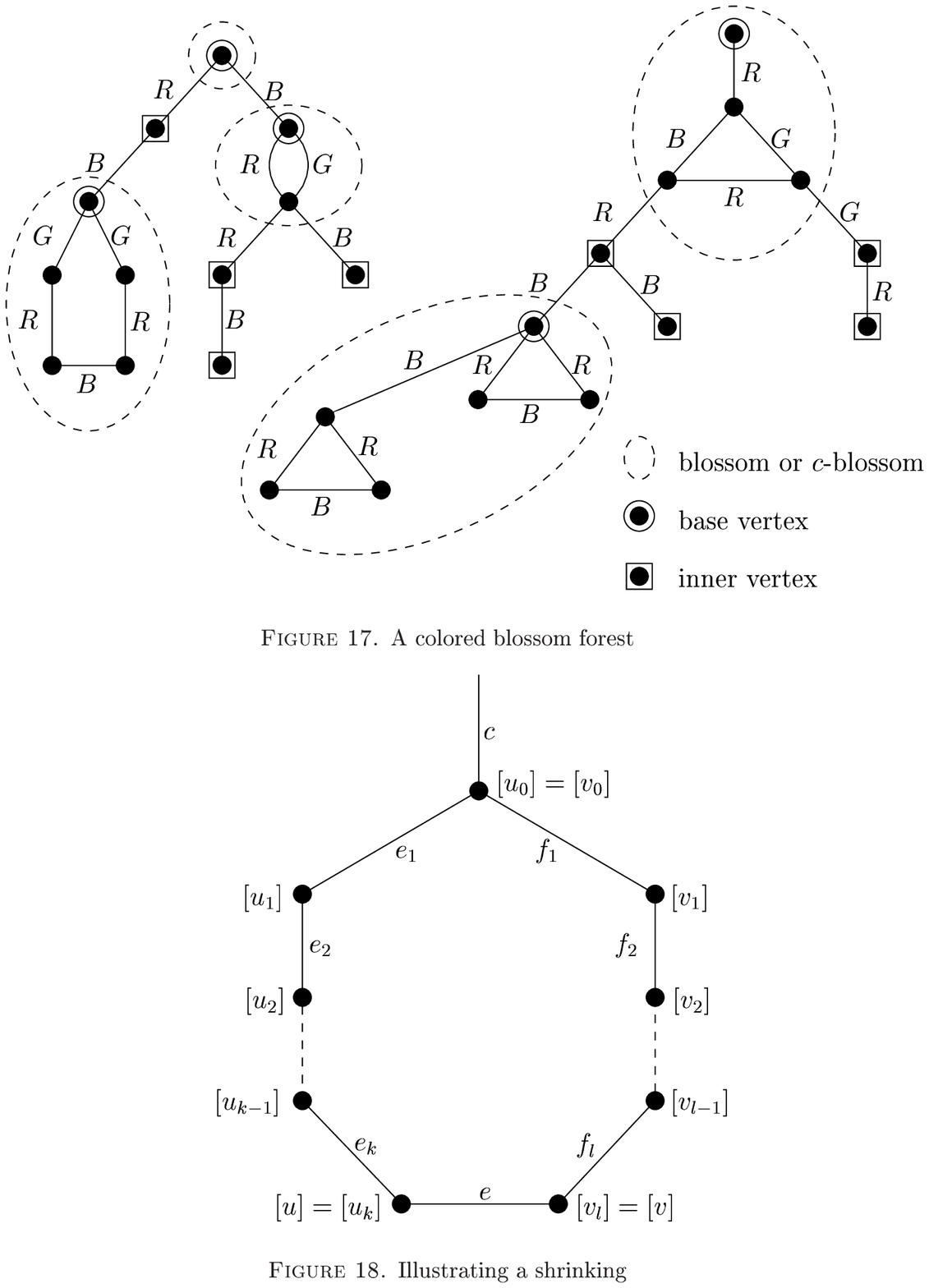}}
   \caption{A colored blossom forest}
   \label{fig_18}
  \end{center}
 \end{figure}

 \bd
 \label{dcbf}
 By a \emph{colored blossom forest} we mean a triple $(I,\pi,F)$,
 where $S \subseteq I \subseteq V$, $\pi$ is a partition of $I$, and
 $F$ is a rooted forest in $G[I]\times \pi$, satisfying the
 following conditions:
 \be
 \item [(i)] $F$ has $\#S$ components and the roots of $F$ are $[s]$, $s \in
 S$;
 \item [(ii)] for each terminal $s \in S$, the induced subgraph $G[[s]]$ is a blossom with
 base $s$;
 \item [(iii)] let $[v]$ be an non root vertex of $F$ satisfying
 $\#[v]=1$,
 let $e$ be the predecessor edge of $[v]$ in $F$, and let
 $e_1,e_2,\ldots ,e_k$ be the edges between $[v]$ and its children
 in $F$; then $e_1,\ldots, e_k$ all have colors different from
 $\cc(e)$;
 \item [(iv)] let $[v]$ be an non root vertex of $F$ satisfying $\#[v] \geq
 2$, let the predecessor edge of $[v]$ in $F$ have color $c$, and let
 its endpoint (in $G$) that is contained in $[v]$ be $u \in [v]$;
 then the induced subgraph $G[[v]]$ is a $c$-blossom with base $u$.
 \ee
 \ed
 Figure~\ref{fig_18} depicts a colored blossom forest, where the
 blocks of $\pi$ are the blossoms, $c$-blossoms, and singleton inner
 vertices indicated in the figure.

 A colored blossom  forest always exists: set $I = S$, take $\pi$ to
 be the trivial partition of $S$ whose blocks are the singletons,
 and take $F$ to be the rooted forest with no edges and roots $[s]$,
 $s \in S$; then $(S,\pi ,F)$ is a colored blossom  forest.

 Given a colored blossom  forest $(I,\pi ,F)$, we classify the
 vertices of $G$  as follows. Vertices in $V - I$ are called
 \emph{out-of-forest vertices}. For $v \in I$, if $\#[v] = 1$ and
 $[v]$ is a non root vertex of $F$, then we call $v$ an \emph{inner
 vertex}. All other vertices in $I$ are called \emph{blossom
 vertices}. If $v$ is an inner or blossom vertex in $V$, we also
 call $[v]$ an inner or blossom vertex in $F$, respectively. In item
 (ii) of Definition~\ref{dcbf}, we say that $s$ is the \emph{base}
 of the blossom vertex $[s]$ of $F$, and in item (iv) of
 Definition~\ref{dcbf}, we say that $u$ is the \emph{base} of the
 blossom vertex $[v]$ of $F$.

 We now define certain positive-length internally alternating trails
 in $G$ w.r.t.\ a colored blossom forest $(I,\pi,F)$.

 Let $u$ be a blossom vertex such that $[u]$ is not a root of $F$.
 Let the predecessor edge of $[u]$ in $F$ have color $c$ and assume
 that $u$ is the base of the blossom vertex $[u]$ of $F$. Pick an
 internally alternating $u$-$u$  trail  of positive length
 in $G[[u]]$ whose first and last edges have colors different from
 $c$ (such a trail is guaranteed by the definition of a colored
 blossom forest), and call it $T(u,F)$.

 Let $u$ be a blossom vertex such that $[u]$ is not a root of $F$.
 Let the predecessor edge of $[u]$ in $F$ have color $c$ and now
 assume that $v$, $v \neq u$, is the base of the blossom vertex
 $[u]$ of $F$. Pick two alternating $u$-$v$ trails in $G[[u]]$ whose first
 edges have different colors and whose last edges have colors
 different from $c$ (such trails are guaranteed by the definition of
 a colored blossom forest), and call them $T_1(u,F)$ and $T_2(u,F)$.

 Let $u \notin S$ be a blossom vertex such that $[u]$ is a root of
 $F$. Let $s \in S$ be the base of the blossom vertex $[u]$ of $F$.
 Pick two alternating $u-s$ trails contained in $G[[s]]$, whose
 first edges have different colors (such trails are guaranteed by
 the definition of a colored blossom forest), and call them
 $T_1(u,F)$ and $T_2(u,F)$.

 The next lemma defines certain internally alternating trails in
 $G$, and states their properties. Some of these trails may be of
 length zero and some may be equal. These conventions prevent some
 case distinctions later on.
 \bl \label{pl} Let $(I,\pi,F)$ be a colored blossom forest, and let
 $u,v$ be vertices of $G$ satisfying the following conditions:
 \begin{itemize}
 \item $u$ and $v$ both belong to $I$;
 \item $[u]$ is a descendent of $[v]$ in $F$;
 \item $v$ is either an inner vertex or the base of the blossom
 vertex $[v]$ of $F$.
 \end{itemize}

 Then $G$ has internally alternating trails $T_1(u,v,F)$,
 $T_2(u,v,F)$ satisfying the following properties.
 \be
 \item [\emph{(i)}] If $u$ is an inner vertex, then $T_1(u,v,F)=T_2(u,v,F)$.
 \item [\emph{(ii)}] The trail $T_1(u,v,F)$ is of length zero precisely when $v=u$.
 The trail $T_2(u,v,F)$ is of length zero precisely when $u=v$ is
 inner or $u=v=s$ for some $s \in S$.
 \item [\emph{(iii)}] The edges of the trails $T_1(u,v,F)$ and $T_2(u,v,F)$ include
 all the edges in the $[u]$-$[v]$ path in $F$.
 \item [\emph{(iv)}] If an edge $e$ in $T_1(u,v,F)$ or $T_2(u,v,F)$ is not in $F$,
 then both its endpoints in $G$ are contained in a blossom vertex of
 $F$ lying on the $[u]$-$[v]$ path in $F$.
 \item [\emph{(v)}] If $u$ is an inner vertex and $v \neq u$, then the first edge of
 $T_1(u,v,F)$ is the predecessor edge of $[u]$ in $F$.
 \item [\emph{(vi)}] If $v$ is an inner vertex and $u \neq v$, then the last edges of
 $T_1(u,v,F)$ and $T_2(u,v,F)$ have colors different from the color
 of the predecessor edge of $[v]$ in $F$.
 \item [\emph{(vii)}] If $u$ is a blossom vertex and $v\neq u$, then the  first
 edges of $T_1(u,v,F),$ and  $T_2(u,v,F)$ have different colors.
 \item [\emph{(viii)}] If $v \notin S$ is the base of a blossom vertex, then the
 first and last edges of $T_2(v,v,F)$ have colors different from the
 color of the predecessor edge of $[v]$ in $F$.
 \item [\emph{(ix)}] If $v \notin S$ is the base of a blossom vertex and $u \neq v$,
 then the last edges of $T_1(u,v,F)$ and $T_2(u,v,F)$ have colors
 different from the color of the predecessor edge of $[v]$ in $F$.
 \ee
 \el
 \pf The proof is by induction on the distance $d([u],[v])$ in $F$
 between the vertices $[u]$ and $[v]$.

 First assume that $d([u],[v]) = 0$, i.e., $u$ and $v$ belong to the
 same block of $\pi$. The following cases arise.
 \newline
 \noi \textbf{Case (a):} $u$ and $v$ are both inner. Define
 $T_1(u,v,F)=T_2(u,v,F)=(u)$, the zero-length trail starting and
 ending at $u$.
 \newline
 \noi \textbf{Case (b):} $u$ and $v$ are both blossom
 vertices and $u \neq  v$. Define
 \[T_1(u,v,F)=T_1(u,F),\quad T_2(u,v,F)=T_2(u,F).\]
 \noi \textbf{Case (c):} $u$ and $v$ are both blossom vertices and $u=v,\;u \notin
 S$. Define
 \[T_1(u,v,F)=(u),\;\;T_2(u,v,F)=T(u,F).\]
 \noi \textbf{Case (d):} $u$ and $v$ are both blossom vertices and $u=v$, $u \in S$. Define
 \[T_1(u,v,F)=T_2(u,v,F)=(u).\]
 It is easily seen that $T_1(u,v,F)$ and $T_2(u,v,F)$ are internally
 alternating and conditions (i)--(ix) in the statement of the lemma
 are satisfied (when restricted to $u,v$ satisfying $d([u],[v])=0$).

 Now assume that $d([u],[v]) > 0$. Then $u \neq v$, and let $e$ be the
 predecessor edge of $[u]$ in $F$. Let the endpoints of $e$ in $G$
 be $x$ and $y$, where $x \in [u]$. Then $d([y],[v])=d([u],[v])-1$ and
 $T_1(y,v,F)$ and $T_2(y,v,F)$ will have been defined. The
 following cases arise.
 \newline
 \noi \textbf{Case (a):} $u$ is inner. We have $u=x$. By induction and (vii) and
 (iii) of Definition~\ref{dcbf}, one of
 \[(u,e,y)*T_1(y,v,F) , \quad (u,e,y)*T_2(y,v,F)\]
 is alternating; define $T_1(u,v,F)=T_2(u,v,F)$ to be that trail
 (breaking ties arbitrarily).
 \newline
 \noi \textbf{Case (b):} $u$ is a blossom vertex  and $u\neq x$. By the definition
 of $T_1(u,F)$, induction and (vii), one of
 \[T_1(u,F)*(x,e,y)*T_1(y,v,F) , \quad T_1(u,F)*(x,e,y)*T_2(y,v,F)\]
 is alternating; define $T_1(u,v,F)$ to be that trail (breaking ties
 arbitrarily). Similarly, one of
 \[T_2(u,F)*(x,e,y)*T_1(y,v,F) , \quad T_2(u,F)*(x,e,y)*T_2(y,v,F)\]
 is alternating; define $T_2(u,v,F)$ to be that trail (breaking
 ties arbitrarily).
 \newline
 \noi \textbf{Case (c):} $u$ is a blossom vertex and $u=x$, $u \notin S$.
 By (iii) of Definition~\ref{dcbf}, induction and (vii), one of
 \[(x,e,y)*T_1(y,v,F) , \quad (x,e,y)*T_2(y,v,F)\]
 is alternating; define $T_1(u,v,F)$ to be that trail (breaking ties
 arbitrarily). Similarly and by the definition of $T(u,F)$, one of
 \[T(u,F)*(x,e,y)*T_1(y,v,F) , \quad T(u,F)*(x,e,y)*T_2(y,v,F)\]
 is alternating; define $T_2(u,v,F)$ to be that trail (breaking
 ties arbitrarily).
 \newline
 \noi \textbf{Case (d):} $u$ is a blossom vertex and $u=x$, $u \in S$. This implies
 $u=v$ and thus $d([u],[v])=0$. So this case cannot occur.

 \noi It is easily checked that $T_1(u,v,F)$ and $T_2(u,v,F)$ are
 alternating and conditions (i)--(ix) in the statement of the lemma
 are satisfied. \epf

 Let $(I,\pi,F)$ be a colored blossom forest and let $e$ be an edge
 in $G$, with endpoints $u$ and $v$ belonging to $I$, that is not an
 edge of $F$. Assume that
 \begin{itemize}
 \item $[u]$ and $[v]$ are not in the same component of $F$.
 \item If $[u]$ is an inner vertex of $F$, then $\cc(e)$ is
 different from the color of the predecessor edge of $[u]$ in $F$.
 \item If $[v]$ is an inner vertex of $F$, then $\cc(e)$ is
 different from the color of the predecessor edge of $[v]$ in $F$.
 \end{itemize}
 In this situation we say that we have a \emph{breakthrough}.

 \bl
 \label{btl}
 Assume we have a breakthrough, with the notation as
 in the preceding paragraph. Let $[s],[t]$, $s,t \in S$, $s \neq t$
 be the roots of the components of $F$ containing $[u]$ and $[v]$,
 respectively. Then $G$ has an alternating $s$-$t$ trail.
 \el
 \pf
 By (vii) of Lemma~\ref{pl} and the definition of breakthrough, for
 some $p,q \in \{1,2\}$,
 \[T_p(u,s,F)^R * (u,e,v) * T_q(v,t,F)\]
 is an alternating $s$-$t$ trail in $G$.
 \epf

 We now discuss three operations on a colored blossom forest
 $(I,\pi,F)$: growing, shrinking and fusing. The operation of fusing
 does not occur in the classical case of searching for augmenting
 paths in a nonbipartite graph.

 Let $e$ be an edge of $G$ between vertices $u$ and $v$ that is not
 an edge of $F$.

 Assume that
 \begin{itemize}
 \item $u \in I$, $v \notin I$.
 \item If $[u]$ is an inner vertex of $F$, then $\cc(e)$ is
 different from the color of the predecessor edge of $[u]$ in $F$.
 \end{itemize}
 Add the singleton block $\{v\}$ to $\pi$ to get a partition $\pi'$
 of $I' = I \cup \{v\}$. Let $F'$ denote the rooted forest in
 $G\times \pi'$ obtained by adding the inner vertex $\{v\}$ to the
 vertices of $F$ and adding the edge $e$ to the set of edges of $F$.
 It is easily seen that $(I',\pi',F')$ is a colored blossom forest.
 We say that $(I',\pi',F')$ is obtained from $(I,\pi,F)$ by
 \emph{growing}.

 We now define the operation of shrinking.

 Let $e$ be an edge of $G$ between vertices $u$ and $v$ that is not
 an edge of $F$.

\begin{figure}
  \begin{center}
  \scalebox{.6}
  {\epsfclipon
  \epsffile[175 200 415 460]{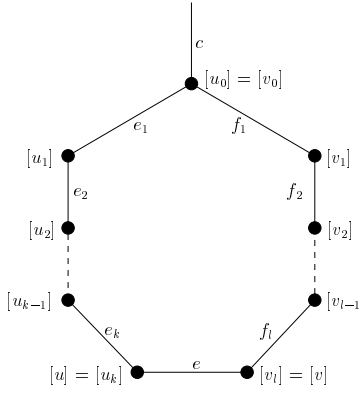}}
   \caption{Illustrating a shrinking}
   \label{fig_19}
  \end{center}
 \end{figure}

 Assume that
 \begin{itemize}
 \item $u, v \in I$.
 \item $[u]$ and  $[v]$ are in the same component of $F$ and
 $[u]\neq [v]$.
 \item If $[u]$ is an inner vertex of $F$, then $\cc(e)$ is
 different from the color of the predecessor edge of $[u]$ in $F$.
 \item If $[v]$ is an inner vertex of $F$, then $\cc(e)$ is
 different from the color of the predecessor edge of $[v]$ in $F$.
 \end{itemize}
 Adding the edge $e$ to $F$ creates a unique cycle $K$.  Let
 $[u_0]=[v_0]$ denote the unique common ancestor in $F$ of $[u]$ and
 $[v]$ that belongs to $K$. Denote the vertices on the $[u_0]$-$[u]$
 path in $F$ by $[u_0],[u_1],\ldots ,[u_{k-1}],[u_k]=[u]$, and
 denote the edge of $F$ between $[u_{i-1}]$ and $[u_i]$ by $e_i$.
 Similarly, denote the vertices on the $[v_0]$-$[v]$ path in $F$ by
 $[v_0],[v_1],\ldots ,[v_{l-1}],[v_l]=[v]$, and denote the edge of
 $F$ between $[v_{i-1}]$ and $[v_i]$ by $f_i$. See
 Figure~\ref{fig_19}. Note that $k$ or $l$ may be zero but $k+l \geq
 1$.

 Replace the blocks $[u_0],\ldots,[u_k],[v_0],\ldots,[v_l]$ of $\pi$
 by their union to obtain a partition $\pi'$ of $I$ with $\pi <
 \pi'$. Define a rooted forest $F'$ in $G \times \pi'$ by throwing
 away the edges $e_1,\ldots,e_k, f_1,\ldots,f_l$ from $F$ (the
 remaining edges in $F$ have obvious endpoints in $G \times \pi'$).
 For convenience, we denote the block of $\pi$ containing $w \in I$
 by $[w]$ and the block of $\pi'$ containing $w$ by $\overline w$.

 We assert that $(I,\pi',F')$ is a colored blossom forest. Clearly,
 conditions~(i) and (iii) in Definition~\ref{dcbf} are satisfied by
 $(I,\pi',F')$, and we need to check condition~(ii) or (iv) only for
 the vertex $\overline u_0$ (according as $[u_0]$ is a root of $F$
 or not). We assume that $[u_0]$ is not a root of $F$ and check
 condition~(iv); the case when $[u_0]$ is a root of $F$ and
 we need to check condition~(ii) is similar, and we omit it.

 Without loss of generality we may assume that if $[u_0]$ is a
 blossom vertex of $F$, then its base is $u_0$. Let $c$ be  the
 color of the predecessor edge of $[u_0]$ in $F$.

 We show that $G[\overline u_0]$ is a $c$-blossom with base $u_0$.
 First we check condition~(ii) of Definition~\ref{dcb}. Let $w \in
 \overline u_0$, $w\neq u_0$. Without loss of generality we may
 assume that $w \in [u_i]$ for some $i=0,\ldots,k$.

 The following two cases arise.
 \newline
 \noi \textbf{Case (a):} $[u_i]$ is a blossom vertex of $F$. By
 (vii) and (ix) of Lemma~\ref{pl}, we have that $T_1(w,u_0,F)$ and
 $T_2(w,u_0,F)$ have first edges of different colors and have last
 edges with colors different from $c$.
 \newline
 \noi \textbf{Case (b):} $[u_i]$ is an inner vertex of $F$. In this
 case $w=u_i$. From the assumption on $[u]$ and $[v]$ and (vii) of
 Lemma~\ref{pl}, it follows that for some $p,q \in \{1,2\}$, the
 $w$-$u_0$ trail
 \[T' = T_p(u,w,F)^R * (u,e,v) * T_q(v,u_0,F)\]
 is alternating. Since $w$ is an inner vertex, it now follows from
 (v) and (vi) of Lemma~\ref{pl} that $T_1(w,u_0,F)$ and $T'$ have
 first edges of different colors, and by (ix) of Lemma~\ref{pl}
 their last edges have colors different from $c$.

 We have verified condition~(ii) in Definition~\ref{dcb}. Now
 consider condition~(iii). As before, for some $p,q \in \{1,2\}$, the
 $u_0$-$u_0$ trail
 \[T_p(u,u_0,F)^R * (u,e,v) * T_q(v,u_0)\]
 of positive length is internally alternating with first and last
 edges of colors different from $c$.

 This completes the proof of the assertion that $(I,\pi',F')$ is a
 colored blossom forest. We say that $(I,\pi',F')$ is obtained from
 $(I,\pi,F)$ by \emph{shrinking}.

 Now we define the operation of fusion. Let $e$ be an edge of $G$
 between vertices $u$ and $v$ that is an edge of $F$. Assume that
 \begin{itemize}
 \item $[u]$ and $[v]$ are (necessarily distinct) blossom vertices of $F$.
 \end{itemize}
 Then $e$ must be the predecessor edge of one of $[u]$ or $[v]$.
 Replace the blocks $[u]$ and $[v]$ of $\pi$ by their union to
 obtain a partition $\pi'$ of $I$. Define a rooted forest $F'$ in $G
 \times \pi'$ by throwing away the edge $e$  from $F$ (the remaining
 edges in $F$ have obvious endpoints in $G \times \pi'$). Again we
 denote the block of $\pi'$ containing a vertex $w \in I$ by
 $\overline w$.

 We assert that $(I,\pi',F')$ is a colored blossom forest. Without
 loss of generality, assume that $e$ is the predecessor edge of
 $[u]$. Clearly, conditions~(i) and (iii) in Definition~\ref{dcbf}
 are satisfied by $(I,\pi',F')$, and we need to check condition~(ii)
 or (iv) only for the vertex $\overline v$ (according as $[v]$ is a
 root of $F$ or not). We assume that $[v]$ is not a root of $F$ and
 check condition~(iv); the case when $[v]$ is a root of $F$ and we
 need to check condition~(ii) is similar, and we omit it.

 Let $x$ be the base of the blossom vertex $[v]$ of $F$, and let $c$
 be the color of the predecessor edge of $[v]$ in $F$. We show that
 $G[\overline x]$ is a $c$-blossom with base $x$. Condition~(iii) in
 Definition~\ref{dcb} clearly holds, since $G[[x]]$ is contained in
 $G[\overline x]$. Condition~(ii) of Definition~\ref{dcb} follows
 from Lemma~\ref{pl}: given that $w \in \overline x$, $w \neq x$,
 the two alternating $w$-$x$ trails $T_1(w,x,F)$ and $T_2(w,x,F)$
 are contained in $G[\overline x]$, have first edges of different
 colors, and have last edges of colors different from $c$.

 This proves the assertion that $(I,\pi',F')$ is a colored blossom
 forest. We say that $(I',\pi',F')$ is obtained from $(I,\pi,F)$ by
 \emph{fusion}.

 We now put a partial order on colored blossom  forests. Given
 colored blossom  forests $\alpha = (I,\pi ,F)$ and $\beta =
 (I',\pi' ,F')$, we say that $\alpha < \beta$ if $I$ is a proper
 subset of $I'$, or $I = I'$ and $\pi < \pi'$ (as partitions, i.e.,
 every block of $\pi$ is contained in a block of $\pi'$ and $\pi
 \neq \pi'$).

 The next theorem is the promised converse of Theorem~\ref{ts}.
 \bt
 \label{cts}
 Suppose that $G$ has no alternating trail connecting distinct
 terminals. Let $(I,\pi ,F)$ be a maximal colored blossom forest.
 Let $A$ be the set of inner vertices of $G$ with respect to $(I,\pi
 ,F)$. For each $c \in C$, define
 \[A(c)=\{u \in A : \text{the predecessor edge of }\, [u] \text{ in } F \text{ has color } c\}.\]
 Then $A$ is a Tutte set with coloring given by  $A=\dot{\cup}_{c \in
 C} A(c)$. \et
 Before proving Theorem~\ref{cts}, we prove the following properties
 of the maximal colored blossom  forest $(I,\pi ,F)$.
 \bl
 \label{cfs}
 Under the assumptions of Theorem~\ref{cts}:
 \be
 \item [\emph{(i)}] No edge in $G$ connects a blossom vertex and an
 out-of-forest vertex.
 \item [\emph{(ii)}] If an edge $e$ in $G$ connects an inner vertex
 and an out-of-forest vertex, then $\cc(e)$ agrees with the color
 of the inner vertex.
 \item [\emph{(iii)}] No edge in $G$ connects blossom vertices
 contained in two different vertices of $F$.
 \item [\emph{(iv)}] If an edge $e$ in $G$ connects two inner vertices,
 then $\cc(e)$ agrees with the color of one of them.
 \item [\emph{(v)}] If an edge $e$ in $G$ connects a blossom vertex
 $v$ and an inner vertex $u$, then $\cc(e)$ agrees with the color
 of $u$, except when $v$ is the base of a blossom and $e$ is the
 predecessor edge of $[v]$ in $F$, in which case the colors are
 different.
 \ee
 \el
 \pf
 If condition~(i) or (ii) does not hold, we can grow the colored
 blossom forest $(I,\pi,F)$, contradicting its maximality.

 We now consider condition~(iii). Let $e$ be an edge in $G$ between
 blossom vertices $u$ and $v$. If $[u]$ and $[v]$ are in different
 components of $F$, then we have a breakthrough which, by
 Lemma~\ref{btl}, contradicts our assumption that there are no
 alternating trails connecting distinct terminals. If $[u]$ and
 $[v]$ are in the same component of $F$ and $e$ is an edge of $F$,
 then we can fuse, contradicting the maximality of $(I,\pi,F)$. If
 $[u]$ and $[v]$ are in the same component of $F$ and $e$ is not an
 edge of $F$, then we can shrink, again contradicting the maximality
 of $(I,\pi,F)$.

 Now we verify condition~(iv). Let $e$ be an edge in $G$ between
 the inner vertices $u$ and $v$. If $e$ is an edge of $F$, then it
 is the predecessor edge of one of $[u]$ and $[v]$, and thus
 $\cc(e)$ agrees with the color of $u$ or $v$. If $e$ is not an edge
 of $F$ and $[u]$ and $[v]$ are in different components of $F$ and
 $\cc(e)$ is different from the colors of $u$ and $v$, then we have
 a breakthrough, a contradiction. If $e$ is not an edge of $F$ and
 $[u]$ and $[v]$ are in the same component of $F$ and $\cc(e)$ is
 different from the colors of $u$ and $v$, then we can shrink, a
 contradiction.

 Finally, consider condition (v). Let $e$ be an edge in $G$
 between a blossom vertex $v$ and an inner vertex $u$. If $e$ is
 an edge of $F$, then either $e$ is the predecessor edge of $[u]$,
 in which case $\cc(e)$ agrees with the color of $u$, or $v$ is the
 base of the blossom vertex $[v]$ of $F$ and $e$ is the predecessor
 edge of $[v]$, in which case $\cc(e)$ is different from the color
 of $u$. If $e$ is not an edge of $F$, then $\cc(e)$ agrees with the
 color of $u$, because otherwise we have a breakthrough (if $[v]$
 and $[u]$ are in different components of $F$) or we can shrink (if
 $[v]$ and $[u]$ are in the same component of $F$).
 \epf

 \pf[Proof of Theorem~\ref{cts}]
 Let $K_1,K_2,\ldots ,K_p$ be the components of the subgraph of $G$
 induced by the out-of-forest vertices. Write $S=\{s_1,s_2, \ldots
 ,s_k\}$ and list the non inner vertices of $F$ as
 \[ [s_1],[s_2],\ldots ,[s_k],[s_{k+1}],\ldots ,[s_l] \]
 (so that the set of blossom vertices in $V$ is precisely $[s_1]
 \cup \cdots \cup [s_l]$). Statements~(i) and (iii) in
 Lemma~\ref{cfs} imply that the components of $G-A$ are precisely
 \[G[[s_1]],\ldots ,G[[s_l]],K_1,\ldots ,K_p.\]
 Condition~(i) in Definition~\ref{dfs} now follows. Statement~(iv)
 in Lemma~\ref{cfs} proves condition~(ii)(c) in
 Definition~\ref{dfs}.

 We verify condition~(ii)(a) in Definition~\ref{dfs} using
 statement~(v) in Lemma~\ref{cfs}, by noting that
 $[s_1],\ldots,[s_k]$ are roots of $F$ and therefore have no
 predecessor edge.

 Condition~(ii)(b) in Definition~\ref{dfs} follows from
 statements~(ii) and (v) in  Lemma~\ref{cfs}: from~(ii) we see that
 there are no mismatched edges with an endpoint in $K_1,\ldots
 ,K_p$, and from~(v) we see that there is exactly one mismatched
 edge with an endpoint in each of $G[[s_{k+1}]],\ldots,G[[s_l]]$.
 \epf

 In the spirit of the Gallai-Edmonds decomposition, the out-of-forest,
 inner, and blossom vertices w.r.t.\ a maximal colored
 blossom  forest can be characterized as follows.

 Define $N(S)$ to be the set of all vertices $t \in V-S$ such that
 for all $s \in S$, there is no alternating $s$-$t$ trail in $G$.

 Define $T(S)$ to be the set of all vertices $t \in V-S$ such that
 for some $s \in S$, there are two alternating $s$-$t$ trails
 whose last edges have different colors.

 For a color $c \in C$, define $I(S,c)$ to be the set of all
 vertices $t \in V-S$ satisfying the following property: there are
 alternating trails starting from $S$ and ending in $t$, and the
 last edges of all such alternating trails have the color $c$. Set
 $I(S)=\dot{\cup}_{c \in C}I(S,c)$.

 Lemma~\ref{os} implies that if $A$ is a Tutte set, then $A
 \subseteq I(S) \cup N(S)$. Theorem~\ref{oib} below and
 Theorems~\ref{ts} and \ref{cts} show that if a Tutte set exists,
 then $I(S)$ is one (indeed, if a Tutte set exists, Theorem~\ref{ts}
 shows that $G$ has no alternating trail connecting distinct
 terminals; then Theorem~\ref{cts} shows that the inner vertices of
 a maximal colored blossom forest form a Tutte set; and
 Theorem~\ref{oib} shows that $I(S)$ are the inner vertices).

 \bt
 \label{oib}
 Suppose that $G$ has no alternating trail connecting distinct
 terminals. Let $(I,\pi ,F)$ be a maximal colored blossom forest.
 Then
 \begin{eqnarray*}
 N(S) & = & \mbox{Out-of-Forest Vertices}
 \\
 I(S,c) & = & \mbox{Inner Vertices with Predecessor Edges having Color $c$}
 \\
 I(S) & = & \mbox{Inner Vertices}
 \\
 S\cup T(S) & = & \mbox{Blossom Vertices}
 \end{eqnarray*}
 \et
 \pf
 Let $X$ be the set of out-of-forest vertices, and let $Y$ be the set of
 blossom vertices. Define $A$ and $A(c)$ as in Theorem~\ref{cts}.
 We first show that $X \subseteq N(S)$, $A(c) \subseteq I(S,c)$, and
 $Y \subseteq
 S\cup T(S)$.

 Write the components of $G-A$ (in the notation of the proof of
 Theorem~\ref{cts}) as
 \[G[[s_1]],\ldots ,G[[s_l]],K_1,\ldots ,K_p.\]
 By Lemma~\ref{cfs}(ii), there is no mismatched edge with an
 endpoint in $K_1,\ldots ,K_p$. It now follows from Theorem~\ref{ts}
 that there is no alternating trail from any vertex in $S$ to any
 out-of-forest vertex, i.e., $X \subseteq N(S)$.

 For each $u \in A(c)$, it follow by Lemmas~\ref{pl}(v) and \ref{os}
 that $u \in I(S,c)$. Thus $A(c) \subseteq I(S,c)$.

 Lemma~\ref{pl}(vii) implies that $Y \subseteq S \cup T(S)$.

 Since $X$, $\dot{\cup}_{c \in C} A(c)$, $Y$ partition $V$ and
 $N(S)$, $I(S)$, $S \cup T(S)$ are disjoint, it follows that
 $X = N(S)$, $A(c) = I(S,c)$, $Y = S \cup T(S)$.
 \epf

 \brm
 The theory presented in this section produces in polynomial
 time either an alternating trail connecting distinct terminals
 or a Tutte set. Indeed, start with the trivial colored blossom forest
 $\pi$ defined after Definition~\ref{dcbf}, and perform growing,
 shrinking and fusing operations in any order until a breakthrough
 or a maximal colored blossom forest is achieved. Discovering
 that one of these operations is possible and performing it or
 discovering a breakthrough takes polynomial time. As for the number
 of operations, initially $\pi$ has $\#S$ blocks. Shrinking and
 fusing decrease the number of blocks of $\pi$ and keep the number
 of out-of-forest vertices constant. Growing increases the first by
 one and decreases the second by one, so at most $\#V - \#S$ growing
 steps can occur in total. It follows that termination must occur
 within $\#V$ operations.
 \erm

 \brm
 We now comment on Tutte's work on the alternating reachability
 problem. Tutte gives a nonalgorithmic solution to a slightly
 different version of the alternating reachability problem. Tutte
 calls the obstructions to the existence of alternating trails
 $r$-barriers \cite[page~331]{t2}. There is a small but important
 difference between our definition of Tutte set and the definition
 of $r$-barrier. If we were to apply Tutte's definition of
 $r$-barrier to the version of alternating reachability considered
 in this paper, then condition (c) in Definition \ref{ts}(ii) would
 read:

 (c) If an edge $e$ connects two vertices of $A$, then these vertices have
 different colors and one of them has the color $\mathcal{C}(e)$,

 instead of our condition paraphrased:

 (c) If an edge $e$ connects two vertices of $A$, then one of them has the color $\mathcal{C}(e)$.

 Thus every $r$-barrier is a Tutte set but not conversely. It is
 easy to find instances of the alternating reachability problem
 where there are no alternating trails connecting distinct
 terminals, but all obstructions are Tutte sets and not
 $r$-barriers. We give such an example for Tutte's original version
 of the alternating reachability problem.

\begin{figure}
 \begin{center}
 \scalebox{.6}
 {\epsfclipon
 \epsffile[280 630 380 750]{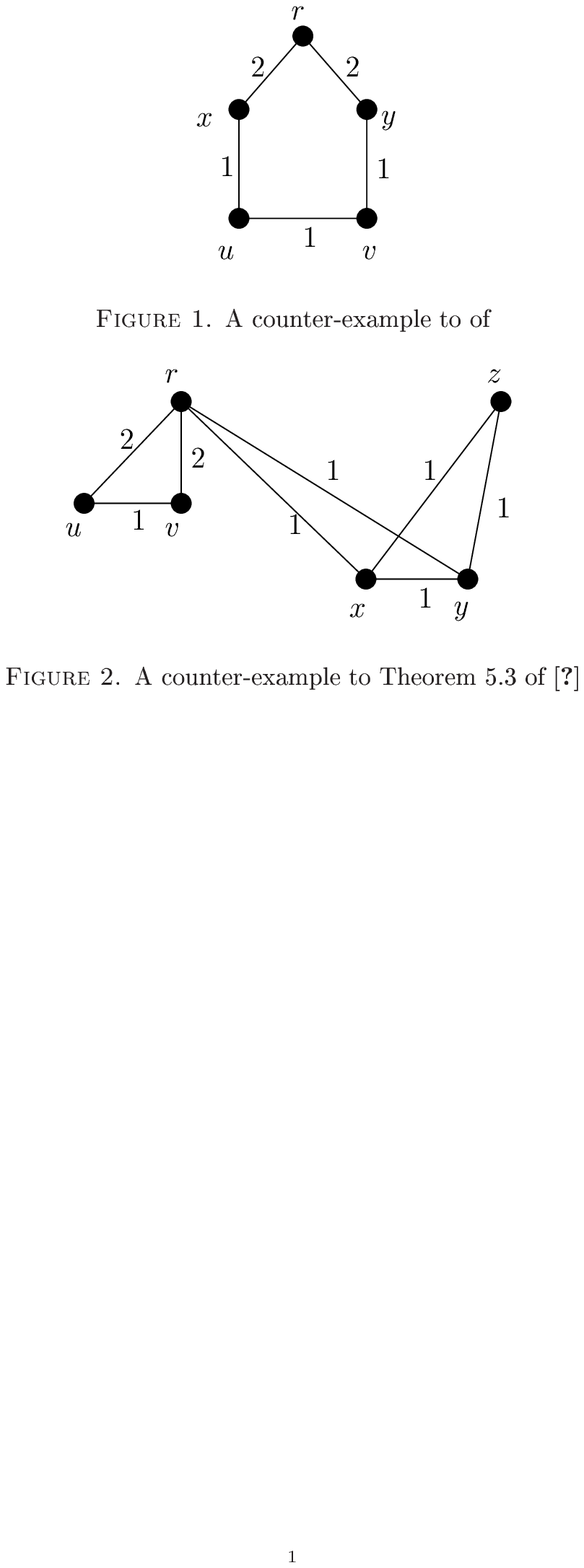}}
  \caption{A counter-example to Theorem~5.1 of \protect\cite{t2}}
  \label{ct_ex_1}
 \end{center}
 \end{figure}

 We follow the definitions and notation of Tutte's paper \cite[page~325--326]{t2}
 without reproducing them here. Consider the
 graph in Figure~\ref{ct_ex_1} (1 and 2 are colors, as in Tutte's
 paper).
 In this graph there are no bicursal edges and the only acursal edge
 is between $u$ and $v$, all vertices are unicursal, and $U_1 = \{
 r,u,v \}$, $U_2 = \{ x,y \}$. This is a counter-example to
 Tutte's Theorem 5.1, since $(U_1,U_2)$ is not an $r$-barrier: the
 edge between $u$ and $v$ violates condition~(ii) in the definition
 of an $r$-barrier.

\begin{figure}
 \begin{center}
 \scalebox{.6}
 {\epsfclipon
 \epsffile[220 490 430 610]{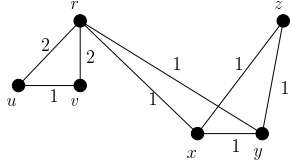}}
  \caption{A counter-example to Theorem~5.3 of \protect\cite{t2}}
  \label{fig_counter_example_tutte}
 \end{center}
 \end{figure}

 Now consider the graph in Figure~\ref{fig_counter_example_tutte}.
 It is easy to convince oneself of the following:
 \be
 \item [(i)] There is no ``coloured path'' in $J(r)$ from $r$ to $z$.
 \item [(ii)] There is no $r$-barrier such that $z$ is a vertex of an
 inaccessible outer component.
 \ee

 In other words, with Tutte's definition of an $r$-barrier, his main result
 Theorem 5.3 is false. The solution is to replace condition (ii)
 in Tutte's definition of an $r$-barrier by the following:

 (ii) If both ends of an edge $A$ of $G$ are inner vertices of $\chi$,
 then one of them has the color of $A$.

 This is precisely our notion of a Tutte set. We believe that the
 only (very minor) error in {\bf \cite{t2}} occurs in the proof of
 Theorem 5.1, and that the theory presented in that paper actually
 yields a Tutte set rather than an $r$-barrier. Tutte has several
 applications of this theorem in his book on graph theory \cite{t1}.
 It is unlikely that these applications need the stronger notion of
 barriers; in most applications, Tutte sets would do. Edges like the
 one between $u$ and $v$ in the first example and between $x$ and
 $y$ in the second are acursal in Tutte's terminology, and therefore
 will not be used by $J(r)$. \erm

 \section{Closed Alternating Trails in Edge-Colored Bridgeless
 Graphs}\label{sec5}
 In this section we give an application of Tutte sets to closed alternating
 trails in edge-colored bridgeless graphs. This result will be used in the
 next section where we determine the inequalities defining the trail
 cone.
 \bt
 \label{catbg}
 Let $G=(V,E)$ be a bridgeless graph with an
 edge-coloring $\cc : E\rar C$. Assume that for each $v \in V$, there
 are two edges with different colors incident at $v$. Then $(G,\cc)$
 has a CAT.
 \et
 \pf Suppose that $(G,\cc)$ has no CAT. Let $e$ be an edge in $G$
 between $s$ and $t$. Form a new edge-colored graph $G'$ from $G$ by
 removing $e$, and adding two new vertices $s'$ and $t'$ and two new
 edges $f$ between $s'$ and $s$ and $g$ between $t'$ and $t$, both
 with the color $\cc(e)$. Since $G$ has no CAT, it follows that $G'$
 has no $s'$-$t'$ alternating trail. Consider the alternating
 reachability problem for $G'$ with $S=\{s',t'\}$ and consider the
 sets $I(S,c)$, $c \in C$ defined in the previous section. Let $A$
 denote the set of inner vertices of a maximal colored blossom
 forest. It follows from Theorem~\ref{oib} that
 $A=I(S)=\dot{\cup}_{c \in C}\;I(S,c)$ is a Tutte set with $s'$ and
 $t'$ in different components of $G' - A$ that are incident with no
 mismatched edges. Since $e$ is not a bridge in $G$, $s'$ and $t'$
 are in the same component of $G'$, and hence $A \neq \emptyset$.

 We first make the following two observations.
 \be
 \item [(i)] There are edges of two different colors incident at every
 vertex of $G'$ other than $s'$ and $t'$.
 \item [(ii)] If $s \notin A$, then $s$ must be in the component of $G'
 - A$ containing $s'$; and similarly if $t \notin A$, then $t$ must
 be in the component of $G' - A$ containing $t'$. \ee

 By definition of $I(S)$, there exist alternating trails in $G'$
 starting at $s'$ or $t'$ and ending at a vertex in $A$. Let $T$
 denote the longest such alternating trail, starting at $s'$ without
 loss of generality and ending at $z \in A$. Let $d$ be the last
 edge of $T$. By Lemma~\ref{os}, $z \in I(S,\cc(d))$. By
 observation~(i) above we can find an edge $h$ of $G'$ incident with
 $z$ satisfying $\cc(d)\neq \cc(h)$. Let $w$ be the other endpoint
 of $h$. We shall now derive a contradiction. The following three
 cases arise.

 \noi \textbf{Case (a):} $T$ contains $h$. Since $\cc(d)\neq
 \cc(h)$, it follows from $z \in I(S,\cc(d))$ that $T$ must have
 traversed $h$ in the direction from $z$ to $w$. Thus, the portion
 of $T$ starting with $h$ and ending with $d$ is a CAT in $G'$.
 Since $s'$ and $t'$ have degree 1 in $G'$, this CAT cannot contain
 the edges $f$ or $g$, so is a CAT in $G$, a contradiction.

 \noi \textbf{Case (b):} $w \in A$ and $T$ does not contain $h$. In
 this case, $T*(z,h,w)$ is an alternating trail that is longer than
 $T$ and ends at a vertex in $A$, a contradiction.

 \noi \textbf{Case (c):} $w \notin A$ and $T$ does not contain $h$.
 Let $H$ be the component of $G' - A$ containing $w$. Since
 $\cc(d)\neq \cc(h)$, it follows that $h$ is mismatched and $H$ does
 not contain $s'$ or $t'$. Hence by observation~(ii) above, $s \in
 A$ or $s \notin H$, and similarly for $t$, so $e$ does not connect
 $A$ and $H$. By the proof of Lemma~\ref{os}, $T$ can enter a
 component of $G' - A$ only via a mismatched edge. It follows (since
 $h$ is the only mismatched edge incident with $H$) that $T$ never
 enters $H$, and so $T$ has no vertices in $H$. Since $z$ is an
 inner vertex and $w$ is a blossom vertex, it follows from
 Lemma~\ref{cfs}(v) that $w$ is the base of a $\cc(h)$-blossom,
 which coincides with $H$ as in the proof of Theorem~\ref{cts}. If
 $h$ were the only edge of $G'$ between $A$ and $H$, then it
 would follow that $h$ is a bridge not only in $G'$ but also in $G$,
 a contradiction. Thus $G$ has some edge $b$, $b\neq h$ between
 $x \in A$ and $y \in H$. Since $T$ has no vertices in $H$, it
 contains neither $h$ nor $b$ nor any edge of $H$. By
 Lemma~\ref{pl}(ii), (vii), (ix) we see that one of
 $T*(z,h,w)*T_1^R(y,w)*(y,b,x)$ and $T*(z,h,w)*T_2^R(y,w)*(y,b,x)$
 is an alternating trail from $s'$ to $A$ that is longer than $T$, a
 contradiction.
 \epf

 \brm
 The algorithm of Section~\ref{sec4} gives a polynomial-time method
 of finding a CAT in an edge-colored graph satisfying the hypothesis
 of Theorem~\ref{catbg}: for each edge $e$ of $G$, construct the
 corresponding graph $G'$ and look for an alternating $s'$-$t'$
 trail. For some edge $e$ we will find an alternating $s'$-$t'$
 trail in $G'$, which easily yields a CAT in $G$. Actually, the
 proof of Theorem~\ref{catbg} gives a method of finding a CAT in $G$
 by solving the alternating reachability problem (with
 $S=\{s',t'\}$) for a single $G'$ (for an arbitrary edge $e$ in
 $G$). We do not give the details.
 \erm

 \brm
 Conjecture~\ref{conj} can be seen as strengthening both the
 hypothesis and conclusion of Theorem~\ref{catbg} (in the case of 2
 colors). The nontrivial part of the conjecture states that if a sum
 of cycles is balanced, then it is a sum of CAT's. We may delete
 bridges and then isolated vertices from $G=(V,E)$, since they are
 irrelevant to the conjecture. Then some positive vector
 $y \in {\mathbb N}^E$ is a sum of cycles. Theorem~\ref{catbg} assumes
 that $E$ (the support of $y$) has edges of both colors at every
 vertex; Conjecture~\ref{conj} assumes more, namely that $y$ is
 balanced. Theorem~\ref{catbg} concludes that $G$ has a CAT, whereas
 Conjecture~\ref{conj} concludes more, namely that $y$ is a sum of
 CAT's.
 \erm

 We conclude this section by deducing the lemma of Giles and Seymour
 on cycles in bridgeless graphs (see \cite{s}) from Theorem
 \ref{catbg}.

 \bt[Giles and Seymour]
 \label{Seymour-Giles}
 Let $G=(V,E)$ be a bridgeless graph and let $\phi : V \rar E$ map
 each vertex $v$ to an edge incident at $v$. Then $G$ has a cycle $C$
 such that for each vertex $w$ of $C$, $\phi(w)$ is an edge of $C$.
 \et
 \pf
 Partition $E$ as $E = E_0 \,\dot{\cup}\, E_1 \,\dot{\cup}\,
 E_2$, where for $i=0,1,2$, $E_i$ consists of those edges $e$ such
 that $e = \phi(v)$ for exactly $i$ endpoints $v$ of $e$.
 Define a 2-colored graph as follows. First subdivide each edge in
 $E_1$ by introducing a new vertex, i.e., for each edge $e \in E_1$
 between $u$ and $w$, introduce a new vertex $v_e$ and replace
 $e$ with the two edges between $u$ to $v_e$ and between $w$ to $v_e$.
 There is an obvious map $\Sigma$ from the edges of the resulting
 graph $G'$ onto $E$ ($\Sigma$ is the identity on $E_0$ and $E_2$
 and takes each subdivided edge onto its parent edge in $E_1$). We
 color the edges in $E_0$ blue and  the edges in $E_2$ red. Now
 consider an edge $e \in E_1$ between $u$ and $w$, say with
 $\phi(u)=e$ and $\phi(w) \neq e$. Then we color the edge between
 $u$ and $v_e$ red and the edge between $w$ and $v_e$ blue.

 It is easily seen that
 \be
 \item [(i)] $G'$ is bridgeless.
 \item [(ii)] Every vertex of $G'$ has positive red and blue
 degrees.
 \item [(iii)] The red edges form a matching in $G'$.
 \ee

 From~(i) and (ii), it follows by Theorem~\ref{catbg} that $G'$ has
 a CAT $T$. From~(iii), $T$ must be an even alternating cycle. Then
 $\Sigma(T)$ is a cycle in $G$ with the required properties.
 \epf

 \section{The Trail Cone}\label{sec6}

 Let $G=(V,E),\; \cc: E \rar \{R,B\}$ be a 2-colored graph. In this
 section we show, using Theorem~\ref{catbg}, that $\tc = \ac \cap
 \zc$. Seymour \cite{s} found the linear inequalities determining
 $\zc$ and our proof is modeled after his. Given a nonempty proper
 subset $X$ of $V$, the subset $D \subseteq E$ of edges between
 $X$ and $V-X$ will be called a \emph{cut}. We say that $X$ and
 $V-X$ are the two \emph{sides} of the cut, and their \emph{sizes}
 are $\#X$ and $\#(V-X)$. Let $D$ be a cut, $e \in D$, and $C$ a
 cycle in $G$. If $C$ contains $e$, then $C$ must also contain an
 edge in $D-\{e\}$. Thus the characteristic vector $\chi(C)$ of $C$
 satisfies the following inequality
 \bq
 \label{cutcondition}
 x(e)  \leq \sum_{f \in D - e}\; x(f),
 \eq
 where we write $D- e$ for $D-\{e\}$. We abbreviate the right-hand
 side of~(\ref{cutcondition}) by $x(D-e)$. We
 call~(\ref{cutcondition}) the \emph{cut condition} for the pair
 $(D,e)$. If it holds with equality, the pair $(D,e)$ is said to be
 \emph{tight for $x$}.

 Seymour \cite{s} proved the following result using
 Theorem~\ref{Seymour-Giles}.
 \bt[Seymour]
 \label{sey}
 $\zc$ is the set of all $x=(x(e): e \in E)$ in ${\mr}^E$ satisfying the
 inequalities
 \begin{alignat}{2}
 \label{cut}
 x(e) & \leq x(D-e), &\qquad& \text{for all cuts }D \text{ and all }e \in D, \\
 \label{nn} x(e) & \geq 0, && \text{for all }e \in E .
 \end{alignat}
 \et

 Vectors satisfying~(\ref{cut})--(\ref{nn}) are said to be
 \emph{cut-admissible for $G$}.

 Given a graph $G=(V,E)$, we let $\kc$ denote the set of all cycles in $G$.
 If $\cc : E\rar\{R,B\}$ is a 2-coloring of $G$, then we denote the set of all CAT's
 in $(G,\cc)$ by $\trc$.
 \bl
 \label{il1}
 Let $G=(V,E)$ be a graph and let $p : E\rar \mq^+$. Let $D$ be a cut
 in $G$, and let $e \in D$ be such that $(D,e)$ is tight for $p$.
 \be
 \item [\emph{(i)}] Suppose $p \in \zc$, which means $p$ can be be expressed as
 \[p = \sum_{C \in\kc}\alpha(C)\chi(C),\qquad \alpha(C) \in \mq^+.\]
 Let $C \in\kc$ with $\alpha(C)>0$. Then $C\cap D$ is either empty or equal
 to $\{e,h\}$ for some $h \in D-e$ (we think of $C$ as a set of edges).
 \item [\emph{(ii)}] Suppose that $\cc:E\rar \{R,B\}$ is a 2-coloring and
 $p \in \tc$, which means $p$ can be be expressed as
 \[p = \sum_{T  \in \trc}\alpha(T)\chi(T),\qquad \alpha(T) \in
 \mq^+.\]
 Let $T \in\trc$ with $\alpha(T)>0$. Then $T\cap D$ is either empty or equal
 to $\{e,h\}$ for some $h \in D-e$.
 \ee
 \el
 \pf
 We prove (i); the proof of (ii) is similar. We have
 \beqn
 \sum_{C \in \kc} \#(C\cap \{e\}) \; \alpha(C) &=& \sum_{C \in \kc,\; e \in C}
 \alpha(C)\\
                                &=& p(e)\\
                                &=& p(D-e)\\
                                &=& \sum_{h \in D-e}\;\;
 \sum_{\stack{C \in\kc}{h \in C}}\alpha(C)\\
 &=& \sum_{C \in \kc}\#(C \cap (D-e)) \; \alpha(C).
 \eeqn
 Since each $C \in \kc$ satisfies $\#(C\cap \{e\}) \leq \#(C \cap
 (D-e))$, it follows that each $C \in \kc$ with $\alpha(C) > 0$
 satisfies $\#(C \cap \{e\}) = \#(C \cap (D-e))$. Since $\#(C \cap
 \{e\}) \in \{0,1\}$, the result follows.
 \epf
 \bl
 \label{il2}
 Let $G=(V,E),\; \cc:E \rar \{R,B\}$ be a 2-colored graph, let $D$
 be a cut in $G$ with sides $X$ and $V-X$, and let $e \in D$ be an
 edge with an endpoint $u_1 \in X$. Let the pair $(D,e)$ be tight
 for the weight function $p : E \rar \mq^+ - \{0\}$.
 \be
 \item [\emph{(i)}] Suppose that $\#X = 1$ and  that $p$ satisfies the balance
 condition at the unique vertex of $X$. Then each edge of $D-e$ has
 color opposite $\cc(e)$. It follows that  for each $T  \in \trc$,
 the intersection $T \cap D$ is either empty or equal to $\{e,h\}$
 for some $h \in D-e$.
 \item [\emph{(ii)}] Suppose that $\#X=2$ and that $p$ satisfies the balance
 condition at both vertices of $X$. Then each edge in $D-e$ with
 color opposite $\cc(e)$ has $u_1$ as an endpoint, and each edge in
 $D-e$ with color $\cc(e)$ does not have $u_1$ as an endpoint. It
 follows that for each $T  \in \trc$, the intersection $T \cap D$ is
 either empty or equal to $\{e,h\}$ for some $h \in D-e$.
 \ee
 \el
 \pf
 (i) This is clear.

 (ii) Let $X=\{u_1,v_1\}$. Without loss of generality we may assume that
 $\cc(e) = R$. Set
 \beqn
 x_1 & = & p(e),
 \\
 x_2 & = & \sum_d \;p(d), \mbox{ where the sum is over all red edges
 $d  \in D-e$ incident with  $u_1$},
 \\
 x_3 & = & \sum_d \;p(d), \mbox{ where the sum is over all blue edges
 $d  \in D$ incident with  $u_1$},
 \\
 x_4 & = & \sum_d \;p(d), \mbox{ where the sum is over all red edges
 $d  \in D$ incident with  $v_1$},
 \\
 x_5 & = & \sum_d \;p(d), \mbox{ where the sum is over all blue edges
 $d  \in D$ incident with  $v_1$},
 \\
 x_6 & = & \sum_d \;p(d), \mbox{ where the sum is over all blue edges
 with both endpoints in $X$},
 \\
 x_7 & = & \sum_d \;p(d), \mbox{ where the sum is over all red edges
 with both endpoints in $X$}.
 \eeqn
 Since the pair $(D,e)$ is tight we have
 \bq
 \label{bal1}
 x_1 = x_2 + x_3 + x_4 + x_5.
 \eq
 The balance condition at $u_1$ gives
 \bq
 \label{bal2}
 x_1+x_2+x_7=x_3+x_6,
 \eq
 and the balance condition at $v_1$ gives
 \bq
 \label{bal3}
 x_5+x_6=x_4+x_7.
 \eq
 Adding~(\ref{bal2}) and (\ref{bal3}) gives
 \[x_1+x_2+x_5=x_3+x_4.\]
 Comparing this equation with~(\ref{bal1}) we get $x_2 + x_5 = 0$,
 and hence $x_2=x_5=0$ by nonnegativity. This implies that there are
 no red edges in $D-e$ incident with $u_1$ and no blue edges in $D$
 incident with $v_1$.
 \epf

 For a graph $G=(V,E)$ and a nonempty proper subset $X$ of $V$, we
 denote by $G_X'$ the graph obtained by shrinking $X$ to a single
 vertex (and deleting the resulting loops).
 \bl
 \label{il3}
 Let $G=(V,E)$ be a graph and let $p:E\rar \mq$. Let $X$ be a
 nonempty proper subset of $V$ and let $p'$ denote $p$ restricted to
 the edges of $G_X'$. If $p$ is cut-admissible for $G$, then $p'$ is
 cut-admissible for $G_X'$.
 \el
 \pf
 This follows since each cut in $G_X'$ is also a cut in $G$.
 \epf

 Let $G=(V,E)$ be a graph with a 2-coloring $\cc : E\rar \{R,B\}$.
 For $X,Y \subseteq V$, we denote by $\nabla_G(X,Y)$ the set of all
 edges of $G$ with one endpoint in $X$ and the other endpoint in
 $Y$. Let $D$ be a cut in $G$ with sides $X$ and $V-X$ of sizes at
 least 3, and let $e \in D$ be an edge between $u_1 \in X$ and
 $u_2 \in V-X$. Let $p : E\rar \mn - \{0\}$ be a weight function.

 Given these data, we define two edge-weighted 2-colored graphs
 $G_X(e)$ (respectively, $G_{V-X}(e)$) by doing the following:
 \begin{itemize}
 \item Delete all edges in $\nabla_G(X,X)$ (respectively, $\nabla_G(V-X,V-X)$).

 \item Replace $X$ (respectively, $V-X$) with $\{u_1,u_1'\}$ (respectively,
 $\{u_2,u_2'\}$), where $u_1',u_2'  \notin V$ are two new vertices.

 \item The endpoints of each edge in $\nabla_G(V-X,V-X) \cup \{e\}$
 (respectively, $\nabla_G(X,X) \cup \{e\}$) remain the same.

 \item The endpoint of each edge $f \in D-e$ in $V-X$ (respectively, $X$) is the same
 as before, and the endpoint in $X$ (respectively, $V-X$) is $u_1$
 (respectively, $u_2$) if $\cc(f)\neq \cc(e)$ and is $u_1'$
 (respectively, $u_2'$) if $\cc(f)=\cc(e)$.

 \item Add a new edge $f_1$ (respectively, $f_2$) between $u_1$ and
 $u_1'$ (respectively, $u_2$ and $u_2'$).

 \item The color of $f_1$ (respectively, $f_2$) is opposite
 $\cc(e)$. All other edges retain their original color.

 \item Define a weight function $p_1$ on the edges of $G_X(e)$ by setting
 $p_1(f_1) = \sum_{h} p(h)$, where the sum is over all $h \in D-e$
 with $\cc(h)=\cc(e)$, and $p_1(h) = p(h)$ for all other edges $h$ of
 $G_X(e)$. Similarly, define a weight function $p_2$ on the edges of
 $G_{V-X}(e)$ by setting $p_2(f_2) = \sum_{h} p(h)$, where the sum
 is over all $h \in D-e$ with $\cc(h)=\cc(e)$, and $p_2(h) = p(h)$ for
 all other edges $h$ of $G_{V-X}(e)$.

 \end{itemize}

 Our restriction on the sizes of $X$ and $V-X$ ensures that $G_X(e)$ and
 $G_{V-X}(e)$ have fewer vertices than $G$.
 \bl
 \label{il4}
 Let $G=(V,E),\;\cc:E\rar \{R,B\}$ be a 2-colored graph, let $D$ be
 a cut in $G$ with sides $X$ and $V-X$ of sizes at least 3, and let
 $e \in D$ be an edge between $u_1 \in X$ and $u_2 \in V-X$. Let
 the pair $(D,e)$ be tight for the weight function $p:E \rar \mn -
 \{0\}$.
 \be
 \item [\emph{(i)}] Suppose that $p$ satisfies the balance condition at
 each vertex of $G$ and is cut-admissible for $G$. Then $p_1$
 (respectively, $p_2$) satisfies the balance condition at each
 vertex of $G_X(e)$ (respectively, $G_{V-X}(e)$) and is
 cut-admissible for $G_X(e)$ (respectively, $G_{V-X}(e)$).
 Furthermore, $(D,e)$ is tight for $p_1$ and $p_2$.

 \item [\emph{(ii)}] Suppose that $p_1$ is a nonnegative integral combination of
 characteristic vectors of CAT's in $G_X(e)$, and $p_2$ is a
 nonnegative integral  combination of characteristic vectors of
 CAT's in $G_{V-X}(e)$. Then $p$ is a nonnegative integral
 combination of characteristic vectors of CAT's in $G$.
 \ee
 \el

 \pf
 (i) From the definition of $p_1(f_1)$ and $p_2(f_2)$ and the
 hypothesis that $(D,e)$ is tight for $p$ and that $p$ satisfies
 the balance condition at each vertex of $G$, it is clear that $p_1$
 and $p_2$ satisfy the balance condition at each vertex of $G_X(e)$
 and $G_{V-X}(e)$, respectively, and that $(D,e)$ is tight for $p_1$ and
 $p_2$. We shall now verify that $p_1$ is cut-admissible for
 $G_X(e)$; the proof for $p_2$ is the same.

 Consider the graph $G_X'$. We retain the name $u_1$ for the vertex
 obtained by shrinking $X$. Let $p'$ be $p$ restricted to the edges
 of $G_X'$. Note that $D$ is a cut in $G_X'$ and $(D,e)$ is tight
 for $p'$. By Lemma~\ref{il3}, $p'$ is cut-admissible for $G_X'$ and
 thus, by Theorem~\ref{sey}, we can write
 \[p' = \sum_{C' \in \kcx} \alpha(C')\chi(C'),\qquad \alpha(C') \in
 \mq^+.\]
 Consider a cycle $C'$ in $\kcx$ with $\alpha(C') > 0$. We obtain a cycle
 $\overline{C'}$ in $G_X(e)$ from $C'$ as follows:
 \begin{itemize}
 \item If $C'$ does not intersect $D$, then $C'$ is a cycle in
 $G_X(e)$ and we set $\overline{C'}=C'$.

 \item If $C'$ intersects $D$ then, by Lemma~\ref{il1}(i), this
 intersection must be  $\{e,h\}$ for some $h \in D-e$. If $u_1$ is
 an endpoint of $h$ in $G_X(e)$ (i.e., if $\cc(e)\neq \cc(h)$) then
 $C'$ is also a cycle in $G_X(e)$ and we set $\overline{C'}=C'$. If
 $u_1$ is not an endpoint of $h$ in $G_X(e)$ (i.e., if
 $\cc(e)=\cc(h)$) we define $\overline{C'}$ to be the cycle in
 $G_X(e)$ obtained from $C'$ by inserting $f_1$ between $e$ and $h$.
 \end{itemize}
 We assert that

 \bq
 \label{assertion}
 p_1 = \sum_{C' \in \kcx} \alpha(C') \chi(\overline{C'}).
 \eq

 To prove the assertion, we examine separately the two kinds of
 edges of $G_X(e)$: edges of $G'_X$ and $f_1$. If $f$ is an edge of
 $G'_X$, then each $C' \in \kcx$ with $\alpha(C') > 0$ satisfies
 $\chi(C')(f) = \chi(\overline{C'})(f)$, hence

 \[\sum_{C' \in \kcx} \alpha(C') \chi(\overline{C'})(f) =
 \sum_{C' \in \kcx} \alpha(C') \chi(C')(f) = p'(f) = p_1(f).\]

 To prove the assertion, it remains to verify

 \[p_1(f_1) = \sum_{C' \in \kcx} \alpha(C') \chi(\overline{C'})(f_1).
 \]

 Consider any $C' \in \kcx$ with $\alpha(C') > 0$ and any $h  \in D -
 e$ with $\cc(h) = \cc(e)$. Then we have, by definition of
 $\overline{C'}$, $\chi(C')(h) = \chi(\overline{C'})(h)  \leq
 \chi(\overline{C'})(f_1)$. It follows that $\chi(C')(h)
 \chi(\overline{C'})(f_1) =  \chi(\overline{C'})(h)$. Thus

 \bq
 \label{hfixed}
 \sum_{C' \in \kcx} \alpha(C') \chi(C')(h) \chi(\overline{C'})(f_1) =
 \sum_{C' \in \kcx} \alpha(C') \chi(\overline{C'})(h).
 \eq

 Now sum~(\ref{hfixed}) over all $h  \in D - e$ such that $\cc(h) =
 \cc(e)$. Since for each $C' \in \kcx$ such that $\alpha(C') > 0$ and
 $f_1  \in \overline{C'}$ there is exactly one edge $h  \in D - e$
 such that $\cc(h) = \cc(e)$ and $h  \in C'$, the left-hand side
 of~(\ref{hfixed}) sums to $\sum_{C' \in \kcx} \alpha(C')
 \chi(\overline{C'})(f_1)$. The right-hand side of~(\ref{hfixed})
 sums to

 \beqn
 \sum_{\stack{h  \in D - e}{\cc(h) = \cc(e)}} \sum_{C' \in \kcx} \alpha(C')
 \chi(\overline{C'})(h)
 &=& \sum_{\stack{h  \in D - e}{\cc(h) = \cc(e)}} \sum_{C' \in \kcx} \alpha(C')
 \chi(C')(h)
 \\
 &=& \sum_{\stack{h  \in D - e}{\cc(h) = \cc(e)}} p'(h)
 \\
 &=& p_1(f_1).
 \eeqn

 This proves the assertion~(\ref{assertion}). Thus $p_1 \in {\mathcal
 Z}(G_X(e))$ and hence $p_1$ is cut-admissible for $G_X(e)$.

 (ii) The hypothesis on $p_1$ (respectively, $p_2$) implies that there is a
 multiset $L_1$ (respectively, $L_2$) of CAT's in $G_X(e)$ (respectively,
 $G_{V-X}(e)$) such that every edge $h$ in $G_X(e)$ (respectively,
 $G_{V-X}(e)$) appears $p_1(h)$ (respectively, $p_2(h)$) times in the
 various CAT's contained in $L_1$ (respectively, $L_2$). We now build a
 multiset $L$ of CAT's in $G$ such that every edge $h$ in $G$ appears
 $p(h)$ times in the CAT's contained in $L$. This will prove the result.

 We begin with some notation. Let $T$ be a CAT in $G$ whose
 intersection with $D$ has exactly 2 edges $e$ and $h$ for some
 $h \in D-e$. Let $h_X$ and $h_{V-X}$ be the endpoints of $h$ in $X$
 and $V-X$, respectively. Then an appropriate cyclic shift of $T$
 must have the form
 \[(u_2,e,u_1)*T_X*(h_X,h,h_{V-X})*T_{V-X},\]
 where $T_X$ is a $u_1$-$h_X$ alternating trail whose vertices are
 in $X$, and $T_{V-X}$ is a $h_{V-X}$-$u_2$ alternating trail whose
 vertices are in $V-X$.

 Consider a CAT in $L_1$ or $L_2$ that intersects $D$. By
 Lemma~\ref{il1}(ii), the intersection of each such CAT with $D$
 must be $\{e,h\}$ for some $h \in D-e$. For $h \in D-e$, let
 $L_1(h)$ (respectively, $L_2(h)$) consist of the CAT's in $L_1$
 (respectively, $L_2$) whose intersection with $D$ is $\{e,h\}$. By
 the definition of $L_1$ and $L_2$, $p_1$ and $p_2$, we have
 $\#L_1(h) = p_1(h) = p(h) = p_2(h) = \#L_2(h)$ for each $h \in
 D-e$. For each $h \in D-e$, fix a bijection $\phi_h:L_1(h)\rar
 L_2(h)$.

 We first take $L$ to be empty and add CAT's to it as follows:
 \begin{itemize}
 \item Add to $L$ all CAT's in $L_1$ whose vertices are contained in $V-X$
 (each such CAT is added the same number of times as it appears in $L_1$).

 \item Add to $L$ all CAT's in $L_2$ whose vertices are contained in $X$.

 \item  For every $h \in D-e$ and every $T \in L_1(h)$ add the CAT
 \[(u_2,e,u_1)*(\phi_h(T))_X*(h_X,h,h_{V-X})*T_{V-X},\]
 to $L$.
 \end{itemize}

 It is easily checked that each edge $h$ in $G$ appears $p(h)$ times
 in the CAT's contained in $L$ (in particular, this holds for $h=e$
 because $(D,e)$ is tight for $p$).
 \epf
 \bt
 Let $G=(V,E),\;\cc:E\rar \{R,B\}$ be a 2-colored graph. Then
 \[\tc = \ac\cap\zc.\]
 \et
 \pf We have already seen in the introduction that $\tc  \subseteq
 \ac\cap\zc$. Consider a nonnegative rational vector
 $q : E\rar{\mq}^+$ that satisfies the balance condition at every
 vertex of $G$ and that is cut-admissible for $G$. We will show
 that $q \in \tc$, which will prove the result. Without loss of
 generality we may assume that $q(e)>0$ for all $e \in E$ (we may
 drop edges $e$ with $q(e)=0$ and maintain the balance condition and
 cut-admissability). The proof is by induction on the pairs
 $(\#V,\#E)$ ordered lexicographically.

 The following two cases arise.

 \textbf{Case (i):} there exist a cut $D$ in $G$ with sides $X$ and
 $V-X$ of sizes at least 3, and an edge $e \in D$ such that $(D,e)$
 is tight for $q$.

 For a suitably large positive integer $N$, the vector $p=Nq$ is
 integral. Thus, by Lemma~\ref{il4}(i), $p_1$ $($respectively,
 $p_2$$)$ satisfies the balance condition at every vertex of
 $G_X(e)$ $($respectively, $G_{V-X}(e)$$)$ and is cut-admissible for
 $G_X(e)$ $($respectively, $G_{V-X}(e)$$)$. Since $G_X(e)$ and
 $G_{V-X}(e)$ have fewer vertices than $G$, we see by induction and
 Lemma~\ref{il4}(ii) that for a suitably large positive integer
 $M$, the vector $Mp$ is a nonnegative integral combination of
 characteristic vectors of CAT's in $G$. It follows that $q \in \tc$,
 as required.

 \textbf{Case (ii):} for each cut $D$ in $G$ with sides $X$ and
 $V-X$ of sizes at least 3 and each $e \in D$, we have $q(e)<q(D-e)$.

 Since $q$ is positive on every edge and cut-admissible for $G$, it
 follows that $G$ is bridgeless. Since $q$ satisfies the balance
 condition at every vertex, it follows from Theorem~\ref{catbg} that
 $(G,\mathcal{C})$ has a CAT $T$.

 Consider the vector $p_t = q - t \chi(T),\;t \geq 0$. Clearly, $p_t$
 is balanced for all $t \geq 0$ and nonnegative for all sufficiently
 small $t > 0$. We assert that there is a positive rational $t_0$
 such that $p_{t_0}$ is cut-admissible for $G$. Indeed, let $D$ be a
 cut in $G$ with sides $X$ and $V-X$, and let $e \in D$. We have the
 following two cases.

 \textbf{Case (a):} $q(e) < q(D-e)$. Clearly $p_t(e) < p_t(D-e)$ for all sufficiently
 small $t > 0$.

 \textbf{Case (b):} $q(e)=q(D-e)$. By assumption, one of $X$ and
 $V-X$, say $X$, has size at most 2. By Lemma~\ref{il2} we see that
 either $T$ contains no edge of $D$ or it contains precisely two
 edges of $D$, $e$ and $h$, for some $h \in D-e$. It follows that
 $p_t(e)=p_t(D-e)$ for all $t \geq 0$.

 From these considerations we see that the maximum value of
 $t$ such that
 \begin{itemize}
 \item $p_t(e) \geq 0$ for all $e \in E$,
 \item $p_t$ is balanced,
 \item $p_t$ is cut-admissible for $G$
 \end{itemize}
 is a positive finite rational $t_0$, as asserted. Set $p=p_{t_0}$.
 The following two subcases arise:

 \textbf{Subcase (ii.1):} $p(f)=0$ for some $f \in E$. By dropping
 $f$ we obtain a graph with the same number of vertices as $G$ but
 with fewer edges, while maintaining balance and cut-admissability.
 Thus by induction $p \in \tc$, and hence $q=p+t_0 \chi(T) \in \tc$,
 as required.

 \textbf{Subcase (ii.2):} $p(f)>0$ for all $f \in E$. From case~(b)
 above we see that $q(e) = q(D-e)$ implies $p(e)=p(D-e)$. Since $p$
 is positive on every edge, it must be that the cutoff determining
 $t_0$ occurs by case~(a) above and not by case~(b) or by the
 requirement that $p_t  \geq 0$. Therefore there is a cut $D^*$ and
 an edge $e^*  \in D^*$ such that $q(e^*) < q(D^*-e^*)$ and
 $p(e^*)=p(D^*-e^*)$. Thus $p$ is a positive rational vector,
 balanced and cut-admissible for $G$, and more pairs $(D,e)$ are
 tight for $p$ than for $q$. We may now repeat the whole argument
 with $p$ in place of $q$. Since the total number of pairs $(D,e)$
 where $D$ a cut in $G$ and $e \in D$ is finite, eventually we will
 reach case (i) or subcase (ii.1).
 \epf

 Finally, we would like to state the following problems.
 In~\cite{BPS1} we saw that the problem of finding an integral
 vector in the intersection of the alternating cone with a box leads
 to the alternating reachability problem. We can ask a similar
 question for the trail cone, but with the integrality restriction
 dropped. Given a 2-colored graph with nonnegative rational upper
 and lower bounds on the edges, is there an augmenting-path-type
 algorithm for either finding a rational vector in the trail cone
 satisfying these bounds, or showing that no such vector exists?

 In {\bf \cite{s}} Seymour makes  the following conjecture for a
 graph $G=(V,E)$: if $y \in (2\mathbb{N})^E \cap {\mathcal Z}(G)$,
 then $y$ is a sum of cycles, i.e, $y$ can be written as a
 nonnegative integer linear combination of characteristic vectors of
 cycles in $G$. A vector in ${\mathcal Z}(G)$ is a nonnegative rational
 combination of characteristic vectors of cycles, i.e., is a
 fractional sum of cycles. So Seymour's conjecture can be stated as
 follows: a fractional sum of cycles that is an even integer on
 every edge is a sum of cycles. Conjecture \ref{conj} states that a
 balanced sum of cycles is a sum of CAT's. Is there any relation
 between these conjectures?

 \end{document}